\documentclass{article}                     

\usepackage{apacite}

\usepackage{graphicx}

\usepackage{booktabs}

\usepackage{float}
\usepackage{multirow}
\usepackage{subfigure}
\usepackage{verbatim}
\usepackage[dvipsnames]{xcolor}
\usepackage{url}
\usepackage{tikz}

\usepackage{algorithmic}
\usepackage[normalem]{ulem}

\usepackage[utf8]{inputenc}
\usepackage{amsmath,amssymb}
\usepackage{amsthm}

\usepackage{pdflscape}
\usepackage{afterpage}
\usepackage{capt-of}
\usepackage{fullpage}

\usepackage[ruled,commentsnumbered,linesnumbered,vlined]{algorithm2e}

\usepackage{indentfirst}
\usepackage{textcomp}

\newcommand{\celsoC}[1]{\textcolor{red}{#1}}

       \makeatletter
\def\@fnsymbol#1{\ensuremath{\ifcase#1\or *\or \mathsection\or \mathparagraph\or \dagger\or \ddagger\or
    \|\or **\or \dagger\dagger
   \or \ddagger\ddagger \else\@ctrerr\fi}}
    \makeatother

\title{Maximum weighted induced forests and trees:\\ New formulations and a computational comparative review}

\author{
Rafael A. Melo {\thanks{Universidade Federal da Bahia, Departamento de Ci\^{e}ncia da Computa\c{c}\~{a}o, Computational Intelligence and Optimization Research Lab (CInO), Salvador, Brazil. ({\tt melo@dcc.ufba.br}) 
}
}
\and Celso C. Ribeiro\thanks{
Universidade Federal Fluminense, Institute of Computing, Niter\'oi, RJ 24210-240, Brazil. ({\tt celso@ic.uff.br})
}
}

\date{\today}

\begin{document}

\maketitle

\begin{abstract}

Given a graph $G=(V,E)$ with a weight $w_v$ associated with each vertex $v\in V$, the maximum weighted induced forest problem (MWIF) consists of encountering a maximum weighted subset $V'\subseteq V$ of the vertices such that $V'$ induces a forest. This NP-hard problem is known to be equivalent to the minimum weighted feedback vertex set problem, which has applicability in a variety of domains. The closely related maximum weighted induced tree problem (MWIT), on the other hand, requires that the subset $V'\subseteq V$ induces a tree.
We propose two new integer programming formulations with an exponential number of constraints and branch-and-cut procedures for MWIF. Computational experiments using benchmark instances are performed comparing several formulations, including the newly proposed approaches and those available in the literature, when solved by a standard commercial mixed integer programming solver. More specifically, five formulations are compared, two compact (i.e., with a polynomial number of variables and constraints) ones and the three others with an exponential number of constraints.
The experiments show that a new formulation for the problem based on directed cutset inequalities for eliminating cycles (DCUT) offers stronger linear relaxation bounds earlier in the search process.
The results also indicate that the other new formulation, denoted tree with cycle elimination (TCYC), outperforms those available in the literature when it comes to the average times for proving optimality for the small instances, especially the more challenging ones.
Additionally, this formulation can achieve much lower average times for solving the larger random instances that can be optimally solved.
Furthermore, we show how the formulations for MWIF can be easily extended for MWIT. Such extension allowed us to analyze the difference between the optimal solution values of the two problems by using some instances of the benchmark set, which is composed of different classes of graphs.
  \\


\end{abstract}

\setlength{\unitlength}{4144sp}%
\begingroup\makeatletter\ifx\SetFigFont\undefined%
\gdef\SetFigFont#1#2#3#4#5{%
  \reset@font\fontsize{#1}{#2pt}%
  \fontfamily{#3}\fontseries{#4}\fontshape{#5}%
  \selectfont}%
\fi\endgroup%

\section{Introduction}
\label{s_intro}

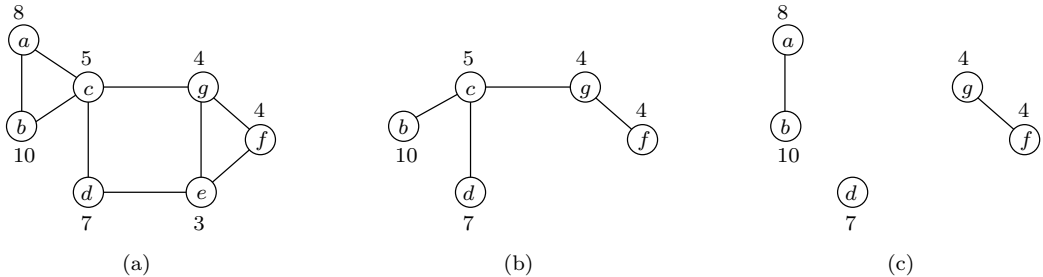
\begin{figure}[ht!]
    \begin{center}
      \subfigure[]{
            \begin{picture}(0,0)%
\includegraphics{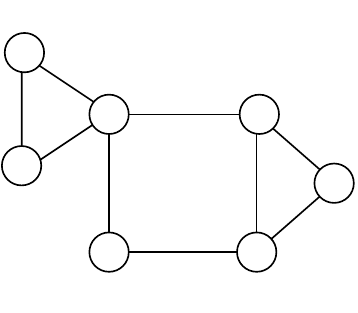}%
\end{picture}%
\setlength{\unitlength}{4144sp}%
\begingroup\makeatletter\ifx\SetFigFont\undefined%
\gdef\SetFigFont#1#2#3#4#5{%
  \reset@font\fontsize{#1}{#2pt}%
  \fontfamily{#3}\fontseries{#4}\fontshape{#5}%
  \selectfont}%
\fi\endgroup%
\begin{picture}(1625,1422)(-92,-1012)
\put(361,-781){\makebox(0,0)[lb]{\smash{{\SetFigFont{8}{9.6}{\rmdefault}{\mddefault}{\updefault}{\color[rgb]{0,0,0}$d$}%
}}}}
\put(-44,-556){\makebox(0,0)[lb]{\smash{{\SetFigFont{8}{9.6}{\rmdefault}{\mddefault}{\updefault}{\color[rgb]{0,0,0}$10$}%
}}}}
\put(361, 29){\makebox(0,0)[lb]{\smash{{\SetFigFont{8}{9.6}{\rmdefault}{\mddefault}{\updefault}{\color[rgb]{0,0,0}$5$}%
}}}}
\put(1396,-286){\makebox(0,0)[lb]{\smash{{\SetFigFont{8}{9.6}{\rmdefault}{\mddefault}{\updefault}{\color[rgb]{0,0,0}$4$}%
}}}}
\put(1036, 29){\makebox(0,0)[lb]{\smash{{\SetFigFont{8}{9.6}{\rmdefault}{\mddefault}{\updefault}{\color[rgb]{0,0,0}$4$}%
}}}}
\put(1036,-961){\makebox(0,0)[lb]{\smash{{\SetFigFont{8}{9.6}{\rmdefault}{\mddefault}{\updefault}{\color[rgb]{0,0,0}$3$}%
}}}}
\put(361,-961){\makebox(0,0)[lb]{\smash{{\SetFigFont{8}{9.6}{\rmdefault}{\mddefault}{\updefault}{\color[rgb]{0,0,0}$7$}%
}}}}
\put(-44,299){\makebox(0,0)[lb]{\smash{{\SetFigFont{8}{9.6}{\rmdefault}{\mddefault}{\updefault}{\color[rgb]{0,0,0}$8$}%
}}}}
\put(-14,131){\makebox(0,0)[lb]{\smash{{\SetFigFont{8}{9.6}{\rmdefault}{\mddefault}{\updefault}{\color[rgb]{0,0,0}$a$}%
}}}}
\put(-24,-387){\makebox(0,0)[lb]{\smash{{\SetFigFont{8}{9.6}{\rmdefault}{\mddefault}{\updefault}{\color[rgb]{0,0,0}$b$}%
}}}}
\put(377,-155){\makebox(0,0)[lb]{\smash{{\SetFigFont{8}{9.6}{\rmdefault}{\mddefault}{\updefault}{\color[rgb]{0,0,0}$c$}%
}}}}
\put(1060,-145){\makebox(0,0)[lb]{\smash{{\SetFigFont{8}{9.6}{\rmdefault}{\mddefault}{\updefault}{\color[rgb]{0,0,0}$g$}%
}}}}
\put(1406,-460){\makebox(0,0)[lb]{\smash{{\SetFigFont{8}{9.6}{\rmdefault}{\mddefault}{\updefault}{\color[rgb]{0,0,0}$f$}%
}}}}
\put(1054,-779){\makebox(0,0)[lb]{\smash{{\SetFigFont{8}{9.6}{\rmdefault}{\mddefault}{\updefault}{\color[rgb]{0,0,0}$e$}%
}}}}
\end{picture}%
            \label{fig:forbCmod32first}
      } \hspace{1cm}
	  \subfigure[]{
            \begin{picture}(0,0)%
\includegraphics{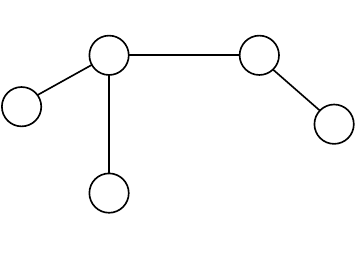}%
\end{picture}%
\setlength{\unitlength}{4144sp}%
\begingroup\makeatletter\ifx\SetFigFont\undefined%
\gdef\SetFigFont#1#2#3#4#5{%
  \reset@font\fontsize{#1}{#2pt}%
  \fontfamily{#3}\fontseries{#4}\fontshape{#5}%
  \selectfont}%
\fi\endgroup%
\begin{picture}(1625,1152)(-92,-1012)
\put(361,-781){\makebox(0,0)[lb]{\smash{{\SetFigFont{8}{9.6}{\rmdefault}{\mddefault}{\updefault}{\color[rgb]{0,0,0}$d$}%
}}}}
\put(-44,-556){\makebox(0,0)[lb]{\smash{{\SetFigFont{8}{9.6}{\rmdefault}{\mddefault}{\updefault}{\color[rgb]{0,0,0}$10$}%
}}}}
\put(361, 29){\makebox(0,0)[lb]{\smash{{\SetFigFont{8}{9.6}{\rmdefault}{\mddefault}{\updefault}{\color[rgb]{0,0,0}$5$}%
}}}}
\put(1396,-286){\makebox(0,0)[lb]{\smash{{\SetFigFont{8}{9.6}{\rmdefault}{\mddefault}{\updefault}{\color[rgb]{0,0,0}$4$}%
}}}}
\put(1036, 29){\makebox(0,0)[lb]{\smash{{\SetFigFont{8}{9.6}{\rmdefault}{\mddefault}{\updefault}{\color[rgb]{0,0,0}$4$}%
}}}}
\put(361,-961){\makebox(0,0)[lb]{\smash{{\SetFigFont{8}{9.6}{\rmdefault}{\mddefault}{\updefault}{\color[rgb]{0,0,0}$7$}%
}}}}
\put(-24,-387){\makebox(0,0)[lb]{\smash{{\SetFigFont{8}{9.6}{\rmdefault}{\mddefault}{\updefault}{\color[rgb]{0,0,0}$b$}%
}}}}
\put(377,-155){\makebox(0,0)[lb]{\smash{{\SetFigFont{8}{9.6}{\rmdefault}{\mddefault}{\updefault}{\color[rgb]{0,0,0}$c$}%
}}}}
\put(1060,-145){\makebox(0,0)[lb]{\smash{{\SetFigFont{8}{9.6}{\rmdefault}{\mddefault}{\updefault}{\color[rgb]{0,0,0}$g$}%
}}}}
\put(1406,-460){\makebox(0,0)[lb]{\smash{{\SetFigFont{8}{9.6}{\rmdefault}{\mddefault}{\updefault}{\color[rgb]{0,0,0}$f$}%
}}}}
\end{picture}%
            \label{fig:forbCmod32second}
	  } \hspace{1cm}
	  \subfigure[]{
            \begin{picture}(0,0)%
\includegraphics{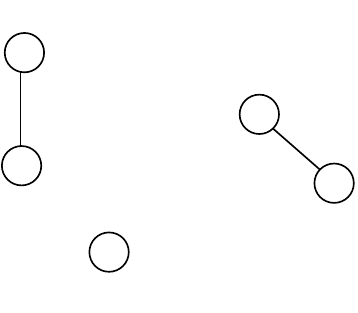}%
\end{picture}%
\setlength{\unitlength}{4144sp}%
\begingroup\makeatletter\ifx\SetFigFont\undefined%
\gdef\SetFigFont#1#2#3#4#5{%
  \reset@font\fontsize{#1}{#2pt}%
  \fontfamily{#3}\fontseries{#4}\fontshape{#5}%
  \selectfont}%
\fi\endgroup%
\begin{picture}(1625,1422)(-92,-1012)
\put(361,-781){\makebox(0,0)[lb]{\smash{{\SetFigFont{8}{9.6}{\rmdefault}{\mddefault}{\updefault}{\color[rgb]{0,0,0}$d$}%
}}}}
\put(-44,-556){\makebox(0,0)[lb]{\smash{{\SetFigFont{8}{9.6}{\rmdefault}{\mddefault}{\updefault}{\color[rgb]{0,0,0}$10$}%
}}}}
\put(1396,-286){\makebox(0,0)[lb]{\smash{{\SetFigFont{8}{9.6}{\rmdefault}{\mddefault}{\updefault}{\color[rgb]{0,0,0}$4$}%
}}}}
\put(1036, 29){\makebox(0,0)[lb]{\smash{{\SetFigFont{8}{9.6}{\rmdefault}{\mddefault}{\updefault}{\color[rgb]{0,0,0}$4$}%
}}}}
\put(361,-961){\makebox(0,0)[lb]{\smash{{\SetFigFont{8}{9.6}{\rmdefault}{\mddefault}{\updefault}{\color[rgb]{0,0,0}$7$}%
}}}}
\put(-44,299){\makebox(0,0)[lb]{\smash{{\SetFigFont{8}{9.6}{\rmdefault}{\mddefault}{\updefault}{\color[rgb]{0,0,0}$8$}%
}}}}
\put(-14,131){\makebox(0,0)[lb]{\smash{{\SetFigFont{8}{9.6}{\rmdefault}{\mddefault}{\updefault}{\color[rgb]{0,0,0}$a$}%
}}}}
\put(-24,-387){\makebox(0,0)[lb]{\smash{{\SetFigFont{8}{9.6}{\rmdefault}{\mddefault}{\updefault}{\color[rgb]{0,0,0}$b$}%
}}}}
\put(1060,-145){\makebox(0,0)[lb]{\smash{{\SetFigFont{8}{9.6}{\rmdefault}{\mddefault}{\updefault}{\color[rgb]{0,0,0}$g$}%
}}}}
\put(1406,-460){\makebox(0,0)[lb]{\smash{{\SetFigFont{8}{9.6}{\rmdefault}{\mddefault}{\updefault}{\color[rgb]{0,0,0}$f$}%
}}}}
\end{picture}%
            \label{fig:forbCmod32third}
	  }
    \end{center}
    \caption{Illustration of maximum weighted induced trees and forests. Subfigure (a) exemplifies an input graph $G$ with seven nodes and $V=\{a,b,c,d,e,f,g\}$. Node weights are represented outside nodes. Subfigure (b) pictures a maximum weighted induced tree $G[V']$, with $V'=\{b,c,d,f,g\}$ and weight 30. Subfigure (c) exemplifies a maximum weighted induced forest $G[V']$, with $V'=\{a,b,d,f,g\}$ and weight 33.}\label{fig:maximumweightedtreesforests}
\end{figure}

Graph theory problems related to encountering induced subgraphs with certain properties have been extensively studied in the literature due to their theoretical interest and various practical applications.
The maximum weighted induced forest problem (MWIF) belongs to this category of problems.
Let us consider a simple and undirected graph $G = (V, E)$ with a set $V$ of vertices and a set $E$ of edges, and a nonnegative weight $w_v$ associated with each vertex $v \in V$. For any subset $V' \subseteq V$ of the vertices, we denote by $G[V']=(V',E')$ the graph induced in $G$ by $V'$, whose edge set $E'$ is formed by all edges in $E$ with both extremities in $V'$, i.e., $E'= \{e = uv \in E \ | \ u,v \in V'\}$. 
The problem thus consists of finding a maximum weighted subset $V'\subseteq V$ inducing a forest $G[V']$, i.e., an acyclic induced subgraph. 
The maximum weighted induced tree problem (MWIT) consists of obtaining a maximum weighted subset $V'\subseteq V$ inducing a tree $G[V']$, i.e., an acyclic and connected induced subgraph. Examples of maximum weighted induced trees and forests are illustrated in Figure~\ref{fig:maximumweightedtreesforests}.

Given a graph $G=(V,E)$, a feedback vertex set is a subset $F\subseteq V$ of its vertices whose removal results in an acycic graph. Note, thus, that $\bar{F} = V \setminus F$ induces a forest.
In this regard, induced forest problems have appeared in several applications in the literature mainly via feedback vertex set problems, which are their equivalent complementary counterparts.  
These include preventing and removing deadlocks in operating systems~\cite{WanLloSof85,CarProSou18}, program verification~\cite{Sha79}, constraint satisfaction~\cite{DecPea88}, study of monopolies in distributed systems~\cite{Pel98}, bayesian inference~\cite{BarGeiNao98}, combinatorial circuit design~\cite{BafBerFuj99}, optical networks~\cite{KleKum01}, parallel systems and distributed computing~\cite{Bos14}, and automated storage/retrieval systems~\cite{GhaXuZha20}.

Several authors have considered theoretical studies and approaches for induced forests and feedback vertex sets. \citeA{Bei97} considered the problem of decycling a graph and introduced the concept of the decycling number of a graph, given by the minimum number of vertices to be removed so that the resulting graph is acyclic (i.e., the size of a minimum cardinality feedback vertex set).
\citeA{BruMafTru00} conducted a polyhedral study and proposed a tabu search metaheuristic for the minimum weighted feedback vertex set problem (MWFVS). 
\citeA{CarCerGenPar05} proposed a linear time algorithm for the MWFVS on a special class of ``diamond'' graphs, as called by the authors.
\citeA{CarCerGen11} implemented a tabu search metaheuristic for the MWFVS which makes use of the algorithms described in~\citeA{CarCerGenPar05} in the local search procedure.
\citeA{CarCerCer14} proposed a memetic algorithm for the MWFVS.
\citeA{ShiXu17} demonstrated graph theoretical results regarding the existence of induced forests with a certain number of vertices.
Very recently, \citeA{MelQueRib21} tackled the MWFVS via the MWIF and proposed compact integer programming formulations and an interated local search (ILS)-based matheuristic for the problem. 
A review of the previous literature on feedback vertex set problems has appeared earlier in~\citeA{FesParRes99}.

To the best of our knowledge, the majority of the works related to encountering induced trees in graphs have concentrated on graph theoretical and computational complexity aspects. \citeA{ErdSakSos86} considered the maximum cardinality induced tree problem and provided bounds for its size which are related to other parameters of the graph.  \citeA{Sco97} proved an open conjecture regarding the existence of induced trees in graphs with large chromatic number. \citeA{Rau07} studied dominating sets inducing trees and provided bounds on the maximum number of vertices of induced trees for certain classes of graphs.
\citeA{DerPic09} considered complexity results for problems related to encountering induced trees covering a prescribed set of vertices. 
\citeA{FoxLohSud09} considered the problem of obtaining large induced trees in $K_r$-free graphs, i.e., those without cliques of size $r$. \citeA{Pfe10} studied large induced trees in triangle-free graphs.
\citeA{ChuSey10} presented a polynomial time algorithm to decide whether there is an induced tree containing three given vertices of a graph.
\citeA{HerMarSch14} proved an upper bound on the difference between the sizes of the maximum cardinality induced forest and the maximum cardinality induced tree in a graph.

Other problems of obtaining connected or induced subgraphs, which are somehow related to the problems considered in our work, were also studied in the recent literature.
\citeA{LjuWeiPfeKlaMutFis06}, \citeA{CosCorLap09}, and \citeA{SieAhmNem20} proposed integer programming approaches for variants of the Steiner
tree problem.
We remark that, although MWIF and MWIT are somehow related to Steiner tree problems, they are very different in nature. This is mainly due to the fact that in the former all vertices are optional.
The maximum weight connected subgraph problem has been studied by several authors~\cite{AlvLjuMut13,RehKoc19,LiuLiJiaHe20}.
\citeA{MarPar16} considered the problem of encountering a maximum balanced induced subgraph.
\citeA{AgrDahHauPin17} tackled the problem of finding maximum $k$-regular induced subgraphs of a graph.
\citeA{MatVerProPas19}, \citeA{BokChiWagWie20}, and~\citeA{MarRib21} considered the problem of encountering the longest induced path.
Besides, integer programming approaches have been successfully applied to several optimization problems related to encountering trees and forests with certain properties~\cite{MelSamUrr16,CarCerCerSil18,CarCerPenRai18}.



The main contributions of this work can be summarized as follows. Firstly, we propose two new formulations for the maximum weighted induced forest problem together with branch-and-cut procedures.
The first is based on constraints that explicitly eliminate cycles in the solution, while the second achieves this goal with the use of the well-known directed cutset inequalities.
To the best of our knowledge, this is the first time these approaches are proposed for MWIF. Besides, although one of the approaches has appeared in the literature for Steiner tree problems, it has not yet been applied and tested for an induced forest problem in which inserting any of the vertices in the solution is optional.
Secondly, we perform extensive computational experiments comparing the proposed formulations with others available in the literature. Thirdly, we show that the proposed formulations can be easily extended to solve MWIT, which, to the best of our knowledge, stand as the first optimization approaches for the problem. The new approaches made it possible to evaluate the difference between the optimal solutions of the two problems for a benchmark set consisting of different classes of graphs.

The remainder of this paper is organized as follows. Section~\ref{sec:formulations} describes the integer programming formulations for the maximum weighted induced forest problem. Section~\ref{sec:branchandcut} details the proposed branch-and-cut approaches.  Section~\ref{sec:experiments} reports the computational experiments. 
Section~\ref{sec:adapting} shows how the new formulations can be extended to deal with the maximum weighted induced tree problem and provides an analysis considering the difference between the optimal solutions of the two problems.
Concluding remarks are discussed in Section
~\ref{sec:concluding}.

\section{Formulations}
\label{sec:formulations}

In this section, we describe three previously existing formulations~\cite{BruMafTru00,MelQueRib21} and
propose two new
formulations for the maximum weighted induced forest problem (MWIF).
The existing formulations are described in Section~\ref{sec:formulationsliterature}, while the new formulations are proposed in Section~\ref{sec:formulationsnew}.

Similarly to \citeA{MelQueRib21}, the input graph for MWIF is slightly transformed such that the problem consists in encountering a tree in this modified graph.
Given the input graph $G=(V,E)$, an alternative transformed graph $G_s=(V_s,E_s)$ is obtained with vertex set $V_s = V \cup \{s\}$ and edge set $E_s = E \cup \{sv \ : \ v \in V \}$.
The weights of all vertices in $G$ remain the same in $G_s$, while the weight of the new vertex is 
$w_s=0$.  
The dummy vertex $s$ is used to join the connected components of the induced forest $G[V']$. 
Some of the formulations consider a directed version of the graph $G_s$, which we denote $\overrightarrow{G}_s = (V_s, A)$. Set $A$ is obtained by creating arcs $(u,v)$ and $(v,u)$ for each edge $uv \in E$ and arcs $(s,v)$ for each vertex $v\in V$. Denote by $A^+(v)$ and $A^-(v)$, respectively, the sets of arcs outgoing from and incoming into vertex $v \in V_s$.
Additionally, throughout the paper, let
$\mathcal{C}$ be the family of all subsets $C\subseteq V$ such that $G[C]$ is a cycle.

\subsection{Existing formulations}
\label{sec:formulationsliterature}
\subsubsection{Cycle elimination formulation}
\label{sec:cycleelimination}

The first described formulation does not have variables corresponding to the edges and is the only one described in this work which does not use the 
transformed graph. It is equivalent to the formulation proposed in \citeA{BruMafTru00} for the minimum weighted feedback vertex set problem, with the only difference being that the variables represent the vertices in the induced forest rather than in the feedback vertex set.  Define the variable $y_v$ to be equal to one if vertex $v\in V$ is in the induced forest $G[V']$, $y_v=0$
otherwise. A cycle elimination-based formulation can be obtained as
\begin{align}
  & \max \sum_{v \in V} w_{v} y_{v}  \label{bru:obj}\\ 
(CYC) \quad & \sum_{v \in C} y_v \leq |C|-1, \ \ \ &  \forall \ C \in \mathcal{C}, \label{bru:01}\\
&  y \in \{0,1\}^{|V|}. \label{bru:02}
\end{align}
The objective function \eqref{bru:obj} maximizes the sum of the weights of the vertices in the forest, i.e., those such that $y_{v} = 1$.
Constraints \eqref{bru:01} ensure the elimination of cycles. Constraints \eqref{bru:02} guarantee the integrality of the variables.

\subsubsection{Flow-based formulation}
\label{sec:flow}

The flow-based formulation~\cite{MelQueRib21}, as well as the two formulations described in Sections~\ref{sec:mtz} and~\ref{sec:directedcutset}, considers the directed graph $\overrightarrow{G}_s = (V_s, A)$.
A flow-based formulation modeling a tree can be obtained by sending one unit of flow from the dummy vertex $s$ to each vertex included in the solution inducing a forest. 
Define variable $f_{uv}$ to represent the flow going from vertex $u$ to vertex $v$, for every arc $(u,v) \in A$.
Let variable $y_{v}$ be defined as before for every $v\in V_s$, and consider variable $x_{uv}$ to be equal to one if arc $(u,v) \in A$ is in the solution, $x_{uv}=0$ otherwise.       
Thus, a flow-based formulation~\cite{MelQueRib21} can be defined as 
\begingroup
\allowdisplaybreaks
\begin{align}
    & \max \sum_{v \in V} w_{v} y_{v} &    \tag{\ref{bru:obj} revisited}\\
(FLOW) \quad & \sum_{(u,v)\in A^-(v)} x_{uv} = y_v, & \forall \ v \in V,  \label{dir:01} \\
& \sum_{(s,v)\in A^+(s)}f_{sv} = \sum_{v\in V} y_v,  \label{dir:02} \\
& \sum_{(u,v) \in A^-(v)} f_{uv} - \sum_{(v,u) \in A^+(v)} f_{vu} = y_{v}, & \forall \ v \in V,  \label{dir:03} \\
& x_{uv} \leq f_{uv} \leq (|V_s| - 1) x_{uv}, & \forall \ (u,v) \in A,  \label{dir:04} \\
& y_s = 1,  \label{dir:05}\\
& x_{uv} + x_{vu}  \leq y_v, & \forall \ (u,v) \in A,   \label{dir:06}\\
& x_{uv} + x_{vu} \geq  y_u + y_v - 1, &   \forall \  uv \in E,  \label{dir:07}\\
& x \in \{0,1\}^{|A|},  \label{dir:08}\\
&  y \in \{0,1\}^{|V_s|}, \label{dir:09}\\
&  f_{uv} \in \mathbb{R}_+^{|A|}. \label{dir:10} 
\end{align}
\endgroup
 Constraints (\ref{dir:01}) ensure there is exactly one arc incoming
 at each selected vertex $v \in V'$. Balance constraints (\ref{dir:02}) and (\ref{dir:03}) are flow conservation constraints for each vertex. Constraints (\ref{dir:04}) link the flow variables $f$ with the $x$ variables. 
 Constraint \eqref{dir:05} forces the dummy node to be in the solution.
 Constraints \eqref{dir:06} and \eqref{dir:07} establish that the solution is an induced subgraph $G[V']$.
  Note that constraints \eqref{dir:06} guarantee that an arc can only be in the solution if both its endpoints are in the solution, while constraints \eqref{dir:07} ensure that if both endpoints of an edge are in the solution, then one of its corresponding arcs must also be in the solution.
  Constraints (\ref{dir:08}),  (\ref{dir:09}), and (\ref{dir:10}) ensure the integrality and nonnegativity of the corresponding variables. 

 As mentioned in \citeA{MelQueRib21}, one could also derive a multicommodity flow-based formulation~\cite{MagWol95} for MWIF, but such type of formulation is usually not viable in practice for large problem instances as even solving its linear relaxation can be challenging~\cite{CosCorLap09,CarCerGauGen13}. 

\subsubsection{Miller-Tucker-Zemlin-based formulation }
\label{sec:mtz}

The Miller-Tucker-Zemlin-based (MTZ-based) formulation for MWIF~\cite{MelQueRib21} also considers the directed graph $\overrightarrow{G}_s$. It adapts the formulation of \citeA{CosCorLap09} which applies the MTZ constraints \cite{MilTucZem60,DesLap91} to solve a variant of the Steiner tree problem.
Define variables $y_v$ and $x_{uv}$ as previously. Additionally, let $\phi_{v}$ be a potential variable associated with each vertex $v \in V_s$. An MTZ-based formulation~\cite{MelQueRib21} can thus be cast as
\begin{align}
  & \max \sum_{v \in V} w_{v} y_{v}  \tag{\ref{bru:obj} revisited}\\  
(MTZ) \quad &  \sum_{(u,v)\in A^-(v)} x_{uv} = y_v, \ \ \ &   \forall \ v \in V,  \tag{\ref{dir:01} revisited} \\
& \phi_{u} - \phi_{v} + |V_s|  x_{uv} + (|V_s|-2) x_{vu} \leq |V_s| - 1, \ \ \ & \forall \ (u,v) \in A,  \label{dir:13}\\
& y_v \leq \phi_{v}  \leq |V_s| - 1, \ \ \  & \forall \  v \in V,  \label{dir:14}\\
& \phi_s = 0, \label{dir:17}\\
& y_s = 1, \tag{\ref{dir:05} revisited}\\
& x_{uv} + x_{vu}  \leq y_v, \ \ \  & \forall \ (u,v) \in A,  \tag{\ref{dir:06} revisited}\\
& x_{uv} + x_{vu} \geq  y_u + y_v - 1, \ \ \  & \forall \ uv \in E,  \tag{\ref{dir:07} revisited}\\
& x \in \{0,1\}^{|A|}, \tag{\ref{dir:08} revisited}\\
&  y \in \{0,1\}^{|V_s|}, \tag{\ref{dir:09} revisited}\\
& \phi \in \mathbb{R}_+^{|V_s|}. \label{dir:19}
\end{align}
Constraints (\ref{dir:13}) guarantee that if an arc $(u,v) \in A$ is in the solution, then vertex $v$ has a potential value $\phi_{v}$ larger than that of $u$, given by $\phi_{u}$.
To see why, note that if $x_{uv} = 1$, then the constraint implies that $\phi_{u} - \phi_{v} \leq - 1$, which is equivalent to $\phi_{v} - \phi_{u} \geq 1$.
Constraints (\ref{dir:14}) determine that if a vertex $v\in V$ is in the solution, then  its potential is neither zero nor greater than $|V_s|-1$. 
Constraint (\ref{dir:17}) sets the potential of the dummy vertex $s$ to zero.  
Constraints (\ref{dir:19}) ensure the nonnegativity of the potential variables.


\subsection{New formulations}
\label{sec:formulationsnew}

\subsubsection{Tree with cycle elimination formulation}
\label{sec:treecycleelimination}

This new undirected formulation considers the modified graph $G_s$ described at the beginning of Section~\ref{sec:formulations}. It guarantees the solution to be acyclic similarly to the cycle elimination formulation described in Section~\ref{sec:cycleelimination}. Nevertheless, it tries to improve the latter from a computational point of view by ensuring that additional constraints are satisfied in an attempt to reduce the computational burden of separating cycle inequalities (which are separated heuristically, as it will be described later in Section~\ref{sec:separationcycle}), as this can potentially allow obtaining better bounds earlier in the search process. Define the variable $y_v$ to be equal to one if vertex $v\in V_s$ is in the induced forest, $y_v=0$ otherwise. Additionally, consider variable $z_e$ to be equal to one if edge $e\in E_s$ is in the induced forest, $z_e=0$ otherwise. The tree with cycle elimination formulation can be obtained as
\begin{align}
  & \max \sum_{v \in V} w_{v} y_{v}  \tag{\ref{bru:obj} revisited}\\
(TCYC) \quad &  \sum_{e \in E_s } z_{e} = \sum_{v \in V} y_v , \ \ \ &    \label{und:01} \\
& \sum_{v \in C} y_v \leq |C|-1, \ \ \ & \forall \ C \in \mathcal{C}, \tag{\ref{bru:01} revisited}\\
& y_s = 1, \tag{\ref{dir:05} revisited} \\
& z_{e} \leq y_v, \ \ \  & \forall \ v \in V,\ e \in \delta(v),  \label{und:04}\\
& z_{e} \geq  y_u + y_v - 1, \ \ \  & \forall \ uv \in E, \label{und:05}\\
& z \in \{0,1\}^{|E_s|}, \label{und:06}\\
&  y \in \{0,1\}^{|V_s|}. \label{und:07}
\end{align}
Constraint \eqref{und:01} ensures the number of edges in the solution equals the number of selected vertices, which represents the number of edges in a tree in $G_s$ including the selected vertices and the dummy vertex. Similarly to constraints~\eqref{dir:06} and \eqref{dir:07}, constraints \eqref{und:04} and \eqref{und:05} force the resulting graph to be induced. 
Constraints \eqref{und:06} and \eqref{und:07} determine the integrality of the $z$ and $y$ variables, correspondingly.

\subsubsection{Directed cutset formulation}
\label{sec:directedcutset}

This new directed cutset formulation also considers the directed graph $\overrightarrow{G}_s$ and guarantees the elimination of cycles using directed cutset connectivity constraints. A feasible solution is an arborescence rooted at the dummy vertex $s$ and, for any partition $\{S,\bar{S}\}$ with $s\in S$ and for any vertex $v\in \bar{S}$, it ensures that if $v$ is in the solution, at least one arc in the cut $(S,\bar{S})$ must be in the solution. Such approach was already applied in the context of Steiner trees \cite{LjuWeiPfeKlaMutFis06,CosCorLap09} but, to the best of our knowledge, has not yet been applied to the problems of obtaining maximum weighted induced forests or trees. Considering variables $y_v$ and $x_{uv}$ as defined in Section~\ref{sec:flow}, a directed cutset formulation can be described as
\begin{align}
 & \max \sum_{v \in V} w_{v} y_{v}  \tag{\ref{bru:obj} revisited}\\
(DCUT) \quad &  \sum_{(u,v)\in A^-(v)} x_{uv} =  y_v , \ \ \ & \forall \ v \in V,   \tag{\ref{dir:01} revisited} \\
& \sum_{(i,j) \in A, i\in S, j\in \bar{S}} x_{ij} \geq y_v, \ \ \ & S \subset V, \ s\in S, \ v \in \bar{S} \label{dircutset:02}\\
& y_s = 1, \tag{\ref{dir:05} revisited}\\
& x_{uv} + x_{vu}  \leq y_v, \ \ \  & \forall \ (u,v) \in A,  \tag{\ref{dir:06} revisited}\\
& x_{uv} + x_{vu} \geq  y_u + y_v - 1, \ \ \  & \forall \ uv \in E,  \tag{\ref{dir:07} revisited}\\
& x \in \{0,1\}^{|A|}, \tag{\ref{dir:08} revisited}\\
&  y \in \{0,1\}^{|V_s|}. \tag{\ref{bru:02} revisited}
\end{align}
Constraints \eqref{dircutset:02} ensure that whenever a vertex $v$ is in the solution, then there is an arc crossing the cut $(S,\bar{S})$, with $s\in S$ and $v\in \bar{S}$.

\newpage

\section{Branch-and-cut approaches}
\label{sec:branchandcut}

This section details the separation procedures for the constraints that are exponential in number.
Furthermore, it describes the clique inequalities, which can be used to strengthen the formulations.



\subsection{Separation of cycle constraints}
\label{sec:separationcycle}

The separation of cycle constraints \eqref{bru:02} takes as input a separation graph $G_{sep}=G[V_{sep}]$ induced by the vertices corresponding to the nonzero $y$ variables in the solution, i.e., $V_{sep} = \{v \in V \ | \ \hat{y}_v > 0\}$. 

The separation for integer solutions is performed exactly using a depth-first search (DFS) algorithm~\cite{Cor09} in $G_{sep}$, which stores every cycle found during the search. After the DFS is finished, all the cycles encountered during its execution are provided to the solver. By representing the graph by an adjacency list, the whole procedure can be implemented to run in time $O(|V|^2 + |E|\times|V|)$.

The separation for fractional solutions is performed heuristically using a DFS algorithm in $G_{sep}$ that considers the vertices in nonincreasing order of the corresponding $\hat{y}$ values. Every cycle found during the search is checked for violation of \eqref{bru:02}, in which case it is stored. After the DFS is finished, all the cycles encountered during the search are given to the solver. By representing the graph by an adjacency list, the whole procedure can be implemented to run in time $O(|V|^2 \lg |V| + |E|\times|V|)$. This complexity comes from the fact that the approach requires sorting the adjacency list in $O(|V|^2 \lg |V|)$ before performing the DFS search with the storage of cycles in $O(|V|^2 + |E|\times|V|)$.

\subsection{Separation of connectivity constraints}

The separation of connectivity constraints \eqref{dircutset:02} also takes as input a separation graph $G_{sep}=G[V_{sep}]$ induced by the vertices corresponding to the nonzero $y$ variables. 

The separation for integer solutions is performed exactly by a breadth-first search (BFS) algorithm~\cite{Cor09} in $G_{sep}$, which initiates at the dummy vertex $s$ and places in the subset $S$ of partition $\{S,\bar{S}\}$ the vertices in the same connected component of $s$. Next, for each vertex in $v\in V_{sep}\setminus S$ an inequality is generated and stored. All the encountered violated inequalities are provided to the solver at the end of the procedure. By representing the graph by an adjacency list, the whole procedure can be implemented to run in $O(|V|^2 + |E|\times|V|)$.

Separation for fractional solutions is performed exactly using maximum flows (minimum cuts) according to the approach described in~\citeA{MagWol95}.
Basically, a new directed graph is built from the solution $(\hat{y},\hat{x})$, where the capacities of the arcs are given by the values of $\hat{x}$. A maximum flow (minimum cut) problem is solved from $s$ to each $v \in V_{sep}$. In case a constraint \eqref{dircutset:02} is found to be violated, it is stored. All the encountered violated inequalities are provided to the solver at the end of the procedure. When using Dinic's algorithm~\cite{Din70} for calculating maximum flows, the procedure can be implemented to run in time $O(|V|^3|E|)$.

\subsection{Separation of clique inequalities}

The clique inequalities for the minimum weighted feedback vertex set problem were proposed in \citeA{BruMafTru00}. They can be defined considering the variables described in our work as follows.
Consider $K_n$ to be any clique, i.e., a complete induced subgraph, with $n\geq 3$. The inequalities 
\begin{equation} \label{ineq:clique}
    \sum_{v \in K_n} y_v \leq 2, \qquad K_n \subset V,
\end{equation}
 are valid for the maximum weighted induced forest problem. 
In fact, more than two vertices of a clique in a solution would induce a cycle.

The separation of clique inequalities \eqref{ineq:clique} is performed heuristically, and also takes as input a separation graph $G_{sep}=G[V_{sep}]$ induced by the vertices corresponding to the nonzero $y$ variables. The heuristic orders the vertices in $V_{sep}$ in nonincreasing order of the corresponding $\hat{y}$ values and in case of ties, in nonincreasing order of degree in $G_{sep}$. The heuristic, then, greedly chooses a vertex which is adjacent to all other vertices already chosen. This procedure is repeated until all vertices are in a maximal clique in $G_{sep}$. Whenever a violated clique is encountered, an attempt to lift it is performed. This is achieved by turning it into a maximal clique in the original graph $G$, whenever possible, via the insertion of vertices which were not in the separation graph. 
This step uses a similar greedy idea to the one used to iteratively build the clique in the separation graph. The difference lies in the fact that it only considers the degree of the vertices as a greedy criterium. By representing the graph by an adjacency list, the whole procedure can be implemented to run in time $O(|V|^2\lg |V| + |E|\times|V|)$.

\section{Computational experiments}
\label{sec:experiments}

In this section, we summarize the computational experiments performed to assess the effectiveness
of the formulations described in Section~\ref{sec:formulations}.
All computational experiments were carried out on a machine running under Ubuntu GNU/Linux, with an Intel(R) Core(TM) i7-4770 CPU @ 3.40GHz processor and 16Gb of RAM. The algorithms were coded in Julia v1.4.2, using JuMP v0.18.6. The formulations were solved using the MIP solver Gurobi 9.0.2.

\subsection{Instances}
\label{sec:instancias}

The experiments were performed using the benchmark instances proposed in \citeA{CarCerGen11} for the maximum weighted feedback vertex set problem, where more details can be encountered. The instances are available in~\citeA{instances}. They consist of random, grid, toroidal, and hypercube graphs, and are classified as small, with $16 \leq |V| \leq 81$, or large,  with $91 \leq |V| \leq 529$. The node weights are uniformly distributed in one of the following intervals: $[10,25]$, $[10,50]$, and $[10,75]$. Similar instances, but with different weights, are organized into instance groups containing five instances each. Each instance group is identified as $X\_n\_m\_low\_up$, where $X$ gives the class of the graph: square grid ($G$), non-square grid ($\mathit{GNQ}$), hypercube ($H$), toroidal ($T$), or random ($R$). 
For random graphs $n$ and $m$ denote, respectively, the numbers of vertices and edges;
for grid and toroidal graphs $n$ and $m$ denote, correspondingly, the number of lines and columns of the grid; 
for hypercube graphs $n=m$ identifies the corresponding $n$-hypercube graph;
the other parameters indicate the lower ($\mathit{low}$) and upper ($\mathit{up}$) bounds on the weights.  
We remark that the results reported throughout this section represent average values over the five instances belonging to the same instance group. 

\subsection{Tested approaches and settings}
\label{sec:testedapproaches}

The following approaches were considered in our experiments:
\begin{itemize}
    \item cycle elimination formulation (CYC) strengthened with clique inequalities; 
    \item compact flow-based formulation (FLOW);
    \item compact MTZ-based formulation (MTZ);
    \item tree with cycle elimination formulation (TCYC) strengthened with clique inequalities; and
    \item directed cutset formulation (DCUT) strengthened with clique inequalities.
\end{itemize}

Note that we follow a common practice in the literature. Namely, FLOW and MTZ are only tested as plain formulations, i.e., without our branch-and-cut procedures. The reason for that is to keep the formulations compact. This makes them a viable option for being easily implemented by a practitioner as one does not have to resort to more advanced algorithmic implementations.

The MIP solver was executed using the standard configurations, except the relative optimality tolerance gap which was set to $10^{-6}$. A time limit of 3600 seconds was imposed for each execution of the MIP solver. All the separation procedures were implemented as callbacks in the MIP solver.
To avoid overloading the formulation with an excessive amount of cuts throughout the search tree, the separations for fractional solutions were only executed at the root node. 
A high quality warm start (i.e., a feasible solution) corresponding to the induced forest generated by the matheuristic MILS$^+$-$mtz$ described in Section 3.6 of~\citeA{MelQueRib21} for the minimum weighted feedback vertex set was offered for each of the executions of the formulations.
We remark that the solver's heuristics were disabled to obtain the linear relaxations at the root node.

We remark that \citeA{BruMafTru00} also showed that certain subset inequalities are valid for the feedback vertex set problem. However, the authors observed in the computational experiments that these inequalities slowed down the solution process and, for this reason, we did not use such inequalities in our experiments.

\subsection{Small instances}
\label{sec:smallinstances}

Tables~\ref{tab:gridinitsol}-\ref{tab:randominitsol} summarize the results for the small instances.
In each of these tables, the first column identifies the instance group. Next, for each of the approaches (CYC, FLOW, MTZ, TCYC, and DCUT), the tables show the average linear relaxation gap achieved by the solver at the end of the root node (glr), given for each instance by $100\times \frac{ub-best}{ub}$, where $ub$ is the bound obtained by the solver at the end of the root node and $best$ is the best solution value encountered by any of the formulations; the number of instances solved to optimality (opt) followed by the average time to solve them (time); and the average relative open gap (gap) in percent considering the instances which were not solved to optimality, given by $100\times \frac{ub-lb}{ub}$, where $lb$ is the best solution value achieved by the formulation.
We note that whenever the optimal solution is known for a given instance, the value glr represents how far (in percent) the linear relaxation bound is from optimality.
Whenever the majority (i.e., at least three) of the formulations are able to solve the same amount of instances, the best average time is shown in bold. Otherwise, the maximum number of instances solved to optimality appears in bold and, in case no instances were solved to optimality by any of the formulations, the smallest gap is represented in bold.

Tables~\ref{tab:gridinitsol}~and~\ref{tab:gridnqinitsol} show that all small grid instances could be solved to optimality using all the five formulations, most of them in low computational times. It can be noted that FLOW achieved the lowest average times for most of the instance groups, while DCUT needed more time on average to solve these instances.
Table~\ref{tab:ipercubeinitsol} shows that, for the small hypercube instances, all formulations but CYC were able to solve all instances to optimality. TCYC achieved the lower average times. 
It can also be noted that the average times for FLOW and MTZ are very similar for these instances.
Besides, considering hypercube instances with the same size, it can be observed that those with weights in the interval $[10,25]$ seem to be more difficult than those with weights in $[10,50]$, which appear to be more challenging than those with weights in $[10,75]$.
Table~\ref{tab:toroinitsol} shows that, for the small toroidal instances, all formulations but CYC solved all instances to optimality, with FLOW achieving lower average times. The most difficult instances in this class appear to be the larger ones with weights varying in the interval $[10,25]$, as it can be observed for instance group $T\_9\_9\_10\_25$.
Table~\ref{tab:randominitsol} shows that for the small random instances, again, all formulations but CYC were able to solve all of them to optimality, with TCYC achieving lower average times, followed by DCUT. It can be seen that, for these instances, TCYC and DCUT outperformed by far MTZ and FLOW, especially for the larger instances in the subset, indicating the benefits of a branch-and-cut approach using the formulations with clique inequalities. Additionally, it can be seen that the density of the instances seems to play an important role in the quality of the bounds provided at the end of the root node, with denser instances showing much higher gaps for all the formulations.

It is noteworthy that DCUT consistently provided lower gaps at the end of the root node and that CYC performed poorly in that respect. For the grid (square and non-square), hypercube and toroidal instances, DCUT was followed, in this order, by FLOW, TCYC, MTZ, and CYC.
For the random instances, this order changed to TCYC, FLOW, MTZ, and CYC. Such change can be explained by the fact that DCUT and TCYC use the branch and cut approaches implementing the separation for clique inequalities.
Besides, it can be noticed that for a few random instance groups, TCYC provided better average bounds, possibly due to the heuristic separation of the clique inequalities.

\begin{landscape}

\begin{table}[H]
\centering
\scriptsize
\caption{Results comparing the formulations for the small square grid instances.}\label{tab:gridinitsol}
\begin{tabular}{l rrrr rrrr rrrr rrrr rrrr}
  \hline
  & \multicolumn{4}{c}{CYC} &  \multicolumn{4}{c}{FLOW} & \multicolumn{4}{c}{MTZ} & \multicolumn{4}{c}{TCYC} & \multicolumn{4}{c}{DCUT} \\ \cmidrule[0.4pt](lr{0.125em}){2-5} \cmidrule[0.4pt](lr{0.125em}){6-9} \cmidrule[0.4pt](lr{0.125em}){10-13}
  \cmidrule[0.4pt](lr{0.125em}){14-17}
  \cmidrule[0.4pt](lr{0.125em}){18-21}
Instance & glr & opt & time & gap & glr & opt & time & gap & glr & opt & time & gap & glr & opt & time & gap & glr & opt & time & gap \\
group & (\%)  &  & (s) & (\%) & (\%)  &  & (s) & (\%) & (\%)  &  & (s) & (\%) & (\%)  &  & (s) & (\%) & (\%)  &  & (s) & (\%) \\ \hline
G\_5\_5\_10\_25  & 3.1 & 5 & 1.4 &    & 0.4 & 5 &  \textbf{0.1}  &    & 2.5 & 5 &  \textbf{0.1}  &    & 1.9 & 5 & 1.1 &    & 0.1 & 5 & 1.6 &   \\ 
  G\_5\_5\_10\_50  & 3.2 & 5 & 1.4 &    & 0.5 & 5 &  \textbf{0.1}  &    & 2.9 & 5 &  \textbf{0.1}  &    & 2.1 & 5 & 1.1 &    & 0.3 & 5 & 1.6 &   \\ 
  G\_5\_5\_10\_75  & 3.0 & 5 & 1.4 &    & 0.4 & 5 &  \textbf{$<$0.1}  &    & 2.8 & 5 & 0.1 &    & 2.0 & 5 & 1.1 &    & 0.2 & 5 & 1.8 &   \\ 
  G\_7\_7\_10\_25  & 7.3 & 5 & 1.5 &    & 0.4 & 5 &  \textbf{0.3}  &    & 1.5 & 5 & 0.6 &    & 1.4 & 5 & 1.2 &    & 0.3 & 5 & 2.1 &   \\ 
  G\_7\_7\_10\_50  & 4.8 & 5 & 1.4 &    & 0.4 & 5 &  \textbf{0.2}  &    & 1.5 & 5 & 0.3 &    & 1.2 & 5 & 1.1 &    & $<$0.1 & 5 & 2.0 &   \\ 
  G\_7\_7\_10\_75  & 6.3 & 5 & 1.5 &    & 0.4 & 5 &  \textbf{0.2}  &    & 1.5 & 5 & 0.4 &    & 1.5 & 5 & 1.2 &    & 0.2 & 5 & 2.1 &   \\ 
  G\_9\_9\_10\_25  & 9.4 & 5 & 5.9 &    & 0.3 & 5 &  \textbf{1.5}  &    & 1.0 & 5 & 11.0 &    & 0.9 & 5 &  \textbf{1.5}  &    & 0.2 & 5 & 3.5 &   \\ 
  G\_9\_9\_10\_50  & 6.3 & 5 & 2.3 &    & 0.6 & 5 & 1.4 &    & 1.3 & 5 & 1.9 &    & 1.2 & 5 &  \textbf{1.3}  &    & 0.3 & 5 & 6.4 &   \\ 
  G\_9\_9\_10\_75  & 6.0 & 5 & 2.0 &    & 0.2 & 5 &  \textbf{0.9}  &    & 0.9 & 5 & 1.4 &    & 1.0 & 5 & 1.2 &    & 0.3 & 5 & 2.8 &   \\ \hline
    Average  & 5.5 &    & 2.1 &    & 0.4 &    & 0.5 &    & 1.8 &    & 1.8 &    & 1.5 &    & 1.2 &    & 0.2 &    & 2.7 &   \\  
  Total &  & 45 &  &  &  & 45 &  &  &  & 45 &  &  &  & 45 &  &  &  & 45 &  &  \\ \hline
\end{tabular}
\end{table}

\begin{table}[H]
\centering
\scriptsize
\caption{Results comparing the formulations for the small non-square grid instances.}\label{tab:gridnqinitsol}
\begin{tabular}{l rrrr rrrr rrrr rrrr rrrr}
  \hline
  & \multicolumn{4}{c}{CYC} &  \multicolumn{4}{c}{FLOW} & \multicolumn{4}{c}{MTZ} & \multicolumn{4}{c}{TCYC} & \multicolumn{4}{c}{DCUT} \\ \cmidrule[0.4pt](lr{0.125em}){2-5} \cmidrule[0.4pt](lr{0.125em}){6-9} \cmidrule[0.4pt](lr{0.125em}){10-13}
  \cmidrule[0.4pt](lr{0.125em}){14-17}
  \cmidrule[0.4pt](lr{0.125em}){18-21}
Instance & glr & opt & time & gap & glr & opt & time & gap & glr & opt & time & gap & glr & opt & time & gap & glr & opt & time & gap \\
group & (\%)  &  & (s) & (\%) & (\%)  &  & (s) & (\%) & (\%)  &  & (s) & (\%) & (\%)  &  & (s) & (\%) & (\%)  &  & (s) & (\%) \\ \hline
GNQ\_8\_3\_10\_25  & 3.1 & 5 & 1.4 &    & 0.6 & 5 &  \textbf{$<$0.1}  &    & 2.7 & 5 & 0.1 &    & 1.6 & 5 & 1.1 &    & 0.1 & 5 & 1.8 &   \\ 
  GNQ\_8\_3\_10\_50  & 1.5 & 5 & 1.4 &    & 0.5 & 5 &  \textbf{0.1}  &    & 2.9 & 5 &  \textbf{0.1}  &    & 1.0 & 5 & 1.1 &    & 0.3 & 5 & 1.8 &   \\ 
  GNQ\_8\_3\_10\_75  & 0.5 & 5 & 1.3 &    & 0.5 & 5 &  \textbf{$<$0.1}  &    & 2.8 & 5 & 0.1 &    & 0.8 & 5 & 1.1 &    & 0.2 & 5 & 1.6 &   \\ 
  GNQ\_9\_6\_10\_25  & 6.5 & 5 & 1.7 &    & 0.2 & 5 &  \textbf{0.2}  &    & 1.2 & 5 & 2.0 &    & 1.1 & 5 & 1.2 &    & 0.1 & 5 & 1.9 &   \\ 
  GNQ\_9\_6\_10\_50  & 5.6 & 5 & 1.5 &    & 0.5 & 5 &  \textbf{0.3}  &    & 1.5 & 5 & 0.7 &    & 1.2 & 5 & 1.2 &    & 0.1 & 5 & 2.1 &   \\ 
  GNQ\_9\_6\_10\_75  & 4.3 & 5 & 1.5 &    & 0.5 & 5 &  \textbf{0.3}  &    & 1.6 & 5 & 0.6 &    & 1.6 & 5 & 1.2 &    & 0.3 & 5 & 2.2 &   \\ 
  GNQ\_12\_6\_10\_25  & 6.1 & 5 & 2.5 &    & 0.3 & 5 &  \textbf{1.2}  &    & 1.1 & 5 & 6.7 &    & 0.9 & 5 & 1.5 &    & 0.3 & 5 & 11.9 &   \\ 
  GNQ\_12\_6\_10\_50  & 5.3 & 5 & 1.8 &    & 0.6 & 5 &  \textbf{1.0}  &    & 1.3 & 5 & 2.2 &    & 1.1 & 5 & 1.2 &    & 0.3 & 5 & 4.0 &   \\ 
  GNQ\_12\_6\_10\_75  & 5.3 & 5 & 1.6 &    & 0.6 & 5 &  \textbf{0.7}  &    & 1.3 & 5 & 1.9 &    & 1.0 & 5 & 1.2 &    & 0.2 & 5 & 3.5 &   \\ \hline
   Average  & 4.3 &    & 1.6 &    & 0.5 &    & 0.4 &    & 1.8 &    & 1.6 &    & 1.1 &    & 1.2 &    & 0.2 &    & 3.4 &   \\ 
  Total &  & 45 &  &  &  & 45 &  &  &  & 45 &  &  &  & 45 &  &  &  & 45 &  &  \\ 
   \hline
\end{tabular}
\end{table}

\begin{table}[H]
\centering
\scriptsize
\caption{Results comparing the formulations for the small hypercube instances.}\label{tab:ipercubeinitsol}
\begin{tabular}{l rrrr rrrr rrrr rrrr rrrr}
  \hline
  & \multicolumn{4}{c}{CYC} &  \multicolumn{4}{c}{FLOW} & \multicolumn{4}{c}{MTZ} & \multicolumn{4}{c}{TCYC} & \multicolumn{4}{c}{DCUT} \\ \cmidrule[0.4pt](lr{0.125em}){2-5} \cmidrule[0.4pt](lr{0.125em}){6-9} \cmidrule[0.4pt](lr{0.125em}){10-13}
  \cmidrule[0.4pt](lr{0.125em}){14-17}
  \cmidrule[0.4pt](lr{0.125em}){18-21}
Instance & glr & opt & time & gap & glr & opt & time & gap & glr & opt & time & gap & glr & opt & time & gap & glr & opt & time & gap \\
group & (\%)  &  & (s) & (\%) & (\%)  &  & (s) & (\%) & (\%)  &  & (s) & (\%) & (\%)  &  & (s) & (\%) & (\%)  &  & (s) & (\%) \\ \hline
H\_4\_10\_25  & 9.7 & 5 & 1.4 &    & 3.1 & 5 &  \textbf{$<$0.1}  &    & 4.3 & 5 & 0.1 &    & 2.4 & 5 & 1.1 &    & 1.5 & 5 & 1.8 &   \\ 
  H\_4\_10\_50  & 5.4 & 5 & 1.4 &    & 2.3 & 5 &  \textbf{$<$0.1}  &    & 4.8 & 5 &  \textbf{$<$0.1}  &    & 2.4 & 5 & 1.1 &    & 0.7 & 5 & 1.8 &   \\ 
  H\_4\_10\_75  & 4.4 & 5 & 1.4 &    & 1.4 & 5 &  \textbf{$<$0.1}  &    & 4.8 & 5 &  \textbf{$<$0.1}  &    & 2.8 & 5 & 1.1 &    & 1.2 & 5 & 1.8 &   \\ 
  H\_5\_10\_25  & 17.2 & 5 & 2.2 &    & 5.8 & 5 & 0.8 &    & 6.8 & 5 &  \textbf{0.7}  &    & 6.9 & 5 & 1.3 &    & 5.4 & 5 & 2.7 &   \\ 
  H\_5\_10\_50  & 11.2 & 5 & 1.5 &    & 2.5 & 5 &  \textbf{0.2}  &    & 3.7 & 5 &  \textbf{0.2}  &    & 4.4 & 5 & 1.2 &    & 2.0 & 5 & 1.9 &   \\ 
  H\_5\_10\_75  & 8.2 & 5 & 1.5 &    & 1.6 & 5 &  \textbf{0.1}  &    & 2.6 & 5 &  \textbf{0.1}  &    & 2.6 & 5 & 1.2 &    & 0.6 & 5 & 1.9 &   \\ 
  H\_6\_10\_25  & 23.2 & 2 & 2296.9 & 2.0 & 5.8 & 5 & 60.7 &    & 6.1 & 5 & 66.3 &    & 6.8 & 5 &  \textbf{22.6}  &    & 5.5 & 5 & 1764.0 &   \\ 
  H\_6\_10\_50  & 15.9 & 5 & 19.6 &    & 3.7 & 5 & 4.8 &    & 3.9 & 5 & 4.1 &    & 4.9 & 5 &  \textbf{2.7}  &    & 3.2 & 5 & 45.8 &   \\ 
  H\_6\_10\_75  & 13.5 & 5 & 7.1 &    & 1.8 & 5 & 1.0 &    & 2.2 & 5 &  \textbf{0.7}  &    & 2.9 & 5 & 1.4 &    & 1.2 & 5 & 3.2 &   \\ \hline
     Average  & 12.1 &    & 259.2 & 2.0 & 3.1 &    & 7.5 &    & 4.4 &    & 8.0 &    & 4.0 &    & 3.7 &    & 2.4 &    & 202.7 &   \\ 
  Total &  & 42 &  &  &  & 45 &  &  &  & 45 &  &  &  & 45 &  &  &  & 45 &  &  \\ 
   \hline
\end{tabular}
\end{table}

\begin{table}[H]
\centering
\scriptsize
\caption{Results comparing the formulations for the small toroidal instances.}\label{tab:toroinitsol}
\begin{tabular}{l rrrr rrrr rrrr rrrr rrrr}
  \hline
  & \multicolumn{4}{c}{CYC} &  \multicolumn{4}{c}{FLOW} & \multicolumn{4}{c}{MTZ} & \multicolumn{4}{c}{TCYC} & \multicolumn{4}{c}{DCUT} \\ \cmidrule[0.4pt](lr{0.125em}){2-5} \cmidrule[0.4pt](lr{0.125em}){6-9} \cmidrule[0.4pt](lr{0.125em}){10-13}
  \cmidrule[0.4pt](lr{0.125em}){14-17}
  \cmidrule[0.4pt](lr{0.125em}){18-21}
Instance & glr & opt & time & gap & glr & opt & time & gap & glr & opt & time & gap & glr & opt & time & gap & glr & opt & time & gap \\
group & (\%)  &  & (s) & (\%) & (\%)  &  & (s) & (\%) & (\%)  &  & (s) & (\%) & (\%)  &  & (s) & (\%) & (\%)  &  & (s) & (\%) \\ \hline
T\_5\_5\_10\_25  & 11.4 & 5 & 1.4 &    & 0.9 & 5 &  \textbf{0.1}  &    & 1.2 & 5 &  \textbf{0.1}  &    & 0.9 & 5 & 1.1 &    & 0.1 & 5 & 1.8 &   \\ 
  T\_5\_5\_10\_50  & 7.0 & 5 & 1.4 &    & 0.6 & 5 &  \textbf{$<$0.1}  &    & 2.1 & 5 & 0.1 &    & 1.9 & 5 & 1.2 &    & 0.2 & 5 & 1.8 &   \\ 
  T\_5\_5\_10\_75  & 7.4 & 5 & 1.4 &    & 0.6 & 5 &  \textbf{$<$0.1}  &    & 2.4 & 5 & 0.1 &    & 2.2 & 5 & 1.1 &    & 0.3 & 5 & 1.8 &   \\ 
  T\_7\_7\_10\_25  & 10.0 & 5 & 4.6 &    & 1.1 & 5 &  \textbf{0.5}  &    & 1.7 & 5 & 0.7 &    & 1.3 & 5 & 1.3 &    & 0.9 & 5 & 2.2 &   \\ 
  T\_7\_7\_10\_50  & 7.4 & 5 & 1.7 &    & 0.8 & 5 &  \textbf{0.3}  &    & 1.7 & 5 & 0.4 &    & 1.7 & 5 & 1.3 &    & 0.6 & 5 & 2.1 &   \\ 
  T\_7\_7\_10\_75  & 5.8 & 5 & 1.6 &    & 0.8 & 5 &  \textbf{0.3}  &    & 1.7 & 5 & 0.4 &    & 1.6 & 5 & 1.2 &    & 0.4 & 5 & 2.2 &   \\ 
  T\_9\_9\_10\_25  & 11.3 & 0 &    & 1.7 & 1.0 & 5 & 5.0 &    & 1.5 & 5 & 17.9 &    & 1.8 & 5 & 10.2 &    & 0.8 & 5 &  \textbf{4.0}  &   \\ 
  T\_9\_9\_10\_50  & 6.3 & 5 & 41.4 &    & 0.6 & 5 &  \textbf{1.0}  &    & 1.1 & 5 & 1.4 &    & 1.1 & 5 & 1.6 &    & 0.4 & 5 & 3.6 &   \\ 
  T\_9\_9\_10\_75  & 6.5 & 5 & 12.7 &    & 0.6 & 5 &  \textbf{1.1}  &    & 1.1 & 5 & 1.2 &    & 1.2 & 5 & 1.5 &    & 0.4 & 5 & 3.6 &   \\   \hline
    Average  & 8.1 &    & 8.3 & 1.7 & 0.8 &    & 0.9 &    & 1.6 &    & 2.5 &    & 1.5 &    & 2.3 &    & 0.4 &    & 2.6 &   \\ 
  Total &  & 40 &  &  &  & 45 &  &  &  & 45 &  &  &  & 45 &  &  &  & 45 &  &  \\ 
   \hline
\end{tabular}
\end{table}

\begin{table}[H]
\centering
\scriptsize
\caption{Results comparing the formulations for the small random instances.}\label{tab:randominitsol}
\begin{tabular}{l rrrr rrrr rrrr rrrr rrrr}
  \hline
  & \multicolumn{4}{c}{CYC} &  \multicolumn{4}{c}{FLOW} & \multicolumn{4}{c}{MTZ} & \multicolumn{4}{c}{TCYC} & \multicolumn{4}{c}{DCUT} \\ \cmidrule[0.4pt](lr{0.125em}){2-5} \cmidrule[0.4pt](lr{0.125em}){6-9} \cmidrule[0.4pt](lr{0.125em}){10-13}
  \cmidrule[0.4pt](lr{0.125em}){14-17}
  \cmidrule[0.4pt](lr{0.125em}){18-21}
Instance & glr & opt & time & gap & glr & opt & time & gap & glr & opt & time & gap & glr & opt & time & gap & glr & opt & time & gap \\
group & (\%)  &  & (s) & (\%) & (\%)  &  & (s) & (\%) & (\%)  &  & (s) & (\%) & (\%)  &  & (s) & (\%) & (\%)  &  & (s) & (\%) \\ \hline
R\_25\_33\_10\_25  & 2.7 & 5 & 1.4 &    & 0.7 & 5 &  \textbf{$<$0.1}  &    & 2.8 & 5 & 0.1 &    & 1.9 & 5 & 1.1 &    & 0.5 & 5 & 1.8 &   \\ 
  R\_25\_33\_10\_50  & 1.4 & 5 & 1.4 &    & 0.5 & 5 &  \textbf{$<$0.1}  &    & 2.6 & 5 & 0.1 &    & 1.1 & 5 & 1.1 &    & 0.2 & 5 & 1.8 &   \\ 
  R\_25\_33\_10\_75  & 1.0 & 5 & 1.4 &    & 0.3 & 5 &  \textbf{$<$0.1}  &    & 2.5 & 5 & 0.1 &    & 0.9 & 5 & 1.1 &    & 0.2 & 5 & 1.4 &   \\ 
  R\_25\_69\_10\_25  & 13.5 & 5 & 1.4 &    & 1.7 & 5 & 0.2 &    & 4.2 & 5 &  \textbf{0.1}  &    & 3.7 & 5 & 1.2 &    & 1.7 & 5 & 1.8 &   \\ 
  R\_25\_69\_10\_50  & 16.5 & 5 & 1.5 &    & 3.8 & 5 & 0.2 &    & 5.7 & 5 &  \textbf{0.1}  &    & 5.1 & 5 & 1.2 &    & 3.4 & 5 & 1.9 &   \\ 
  R\_25\_69\_10\_75  & 14.7 & 5 & 1.5 &    & 2.5 & 5 &  \textbf{0.1}  &    & 4.5 & 5 &  \textbf{0.1}  &    & 4.3 & 5 & 1.2 &    & 2.0 & 5 & 1.9 &   \\ 
  R\_25\_204\_10\_25  & 31.9 & 5 & 1.7 &    & 26.4 & 5 & 1.7 &    & 28.3 & 5 &  \textbf{0.9}  &    & 21.8 & 5 & 1.5 &    & 19.9 & 5 & 2.3 &   \\ 
  R\_25\_204\_10\_50  & 31.9 & 5 & 1.6 &    & 29.2 & 5 & 1.3 &    & 30.7 & 5 &  \textbf{0.7}  &    & 21.0 & 5 & 1.4 &    & 15.4 & 5 & 2.3 &   \\ 
  R\_25\_204\_10\_75  & 30.2 & 5 & 1.6 &    & 27.7 & 5 & 1.3 &    & 29.5 & 5 &  \textbf{0.7}  &    & 19.4 & 5 & 1.4 &    & 15.4 & 5 & 2.2 &   \\ 
  R\_50\_85\_10\_25  & 8.4 & 5 & 2.1 &    & 0.3 & 5 &  \textbf{0.1}  &    & 1.3 & 5 & 0.4 &    & 1.3 & 5 & 1.2 &    & 0.1 & 5 & 1.9 &   \\ 
  R\_50\_85\_10\_50  & 7.7 & 5 & 1.8 &    & 0.7 & 5 &  \textbf{0.2}  &    & 1.5 & 5 & 0.3 &    & 1.6 & 5 & 1.2 &    & 0.4 & 5 & 2.1 &   \\ 
  R\_50\_85\_10\_75  & 7.9 & 5 & 1.6 &    & 0.5 & 5 &  \textbf{0.2}  &    & 1.5 & 5 & 0.3 &    & 1.4 & 5 & 1.2 &    & 0.2 & 5 & 2.0 &   \\ 
  R\_50\_232\_10\_25  & 39.8 & 5 & 77.5 &    & 6.1 & 5 & 5.3 &    & 7.2 & 5 & 3.0 &    & 5.7 & 5 &  \textbf{2.0}  &    & 6.1 & 5 & 4.0 &   \\ 
  R\_50\_232\_10\_50  & 39.4 & 5 & 28.7 &    & 6.4 & 5 & 4.7 &    & 7.4 & 5 & 2.2 &    & 6.0 & 5 &  \textbf{1.7}  &    & 6.0 & 5 & 4.2 &   \\ 
  R\_50\_232\_10\_75  & 40.8 & 5 & 27.5 &    & 5.8 & 5 & 4.1 &    & 6.9 & 5 & 2.2 &    & 5.4 & 5 &  \textbf{1.7}  &    & 5.8 & 5 & 4.1 &   \\ 
  R\_50\_784\_10\_25  & 72.6 & 5 & 231.0 &    & 46.4 & 5 & 53.8 &    & 46.6 & 5 & 32.3 &    & 40.3 & 5 &  \textbf{8.4}  &    & 27.1 & 5 & 9.8 &   \\ 
  R\_50\_784\_10\_50  & 73.0 & 5 & 172.1 &    & 47.4 & 5 & 32.6 &    & 48.1 & 5 & 15.9 &    & 39.7 & 5 &  \textbf{8.5}  &    & 27.1 & 5 & 9.0 &   \\ 
  R\_50\_784\_10\_75  & 72.0 & 5 & 133.3 &    & 45.8 & 5 & 31.2 &    & 45.5 & 5 & 13.8 &    & 37.6 & 5 &  \textbf{7.8}  &    & 25.2 & 5 & 8.2 &   \\ 
  R\_75\_157\_10\_25  & 13.9 & 5 & 804.8 &    & 0.9 & 5 &  \textbf{1.0}  &    & 1.4 & 5 & 2.8 &    & 1.8 & 5 & 3.1 &    & 0.8 & 5 & 3.0 &   \\ 
  R\_75\_157\_10\_50  & 15.7 & 5 & 195.9 &    & 1.1 & 5 &  \textbf{0.9}  &    & 1.7 & 5 & 1.1 &    & 2.0 & 5 & 2.2 &    & 0.8 & 5 & 3.4 &   \\ 
  R\_75\_157\_10\_75  & 16.6 & 5 & 62.8 &    & 1.1 & 5 & 1.2 &    & 1.7 & 5 &  \textbf{0.9}  &    & 2.0 & 5 & 1.7 &    & 0.9 & 5 & 5.8 &   \\ 
  R\_75\_490\_10\_25  & 48.5 & 0 &    & 13.7 & 9.7 & 5 & 129.1 &    & 10.4 & 5 & 69.2 &    & 8.0 & 5 &  \textbf{12.5}  &    & 10.1 & 5 & 41.4 &   \\ 
  R\_75\_490\_10\_50  & 50.6 & 0 &    & 15.6 & 12.8 & 5 & 167.6 &    & 12.5 & 5 & 70.7 &    & 10.8 & 5 &  \textbf{20.1}  &    & 12.8 & 5 & 56.7 &   \\ 
  R\_75\_490\_10\_75  & 49.7 & 0 &    & 9.8 & 11.3 & 5 & 81.5 &    & 11.5 & 5 & 32.2 &    & 9.2 & 5 &  \textbf{6.3}  &    & 11.3 & 5 & 35.4 &   \\ 
  R\_75\_1739\_10\_25  & 79.0 & 0 &    & 47.2 & 58.2 & 5 & 808.2 &    & 58.8 & 5 & 372.0 &    & 56.4 & 5 &  \textbf{133.0}  &    & 36.9 & 5 & 138.7 &   \\ 
  R\_75\_1739\_10\_50  & 78.6 & 0 &    & 38.3 & 58.3 & 5 & 690.5 &    & 58.3 & 5 & 331.5 &    & 53.2 & 5 &  \textbf{100.0}  &    & 35.1 & 5 & 103.4 &   \\ 
  R\_75\_1739\_10\_75  & 78.2 & 0 &    & 36.2 & 57.0 & 5 & 695.4 &    & 57.3 & 5 & 298.6 &    & 47.4 & 5 & 97.7 &    & 33.3 & 5 &  \textbf{91.5}  &   \\ \hline
       Average  & 34.7 &    & 83.5 & 26.8 & 17.1 &    & 100.5 &    & 18.2 &    & 46.4 &    & 15.2 &    & 15.6 &    & 11.1 &    & 20.1 &   \\  
  Total &  & 105 &  &  &  & 135 &  &  &  & 135 &  &  &  & 135 &  &  &  & 135 &  &  \\ 
   \hline
\end{tabular}
\end{table}

\end{landscape}

\begin{figure}[ht!]
\begin{center}
     \includegraphics[width=0.8\textwidth]{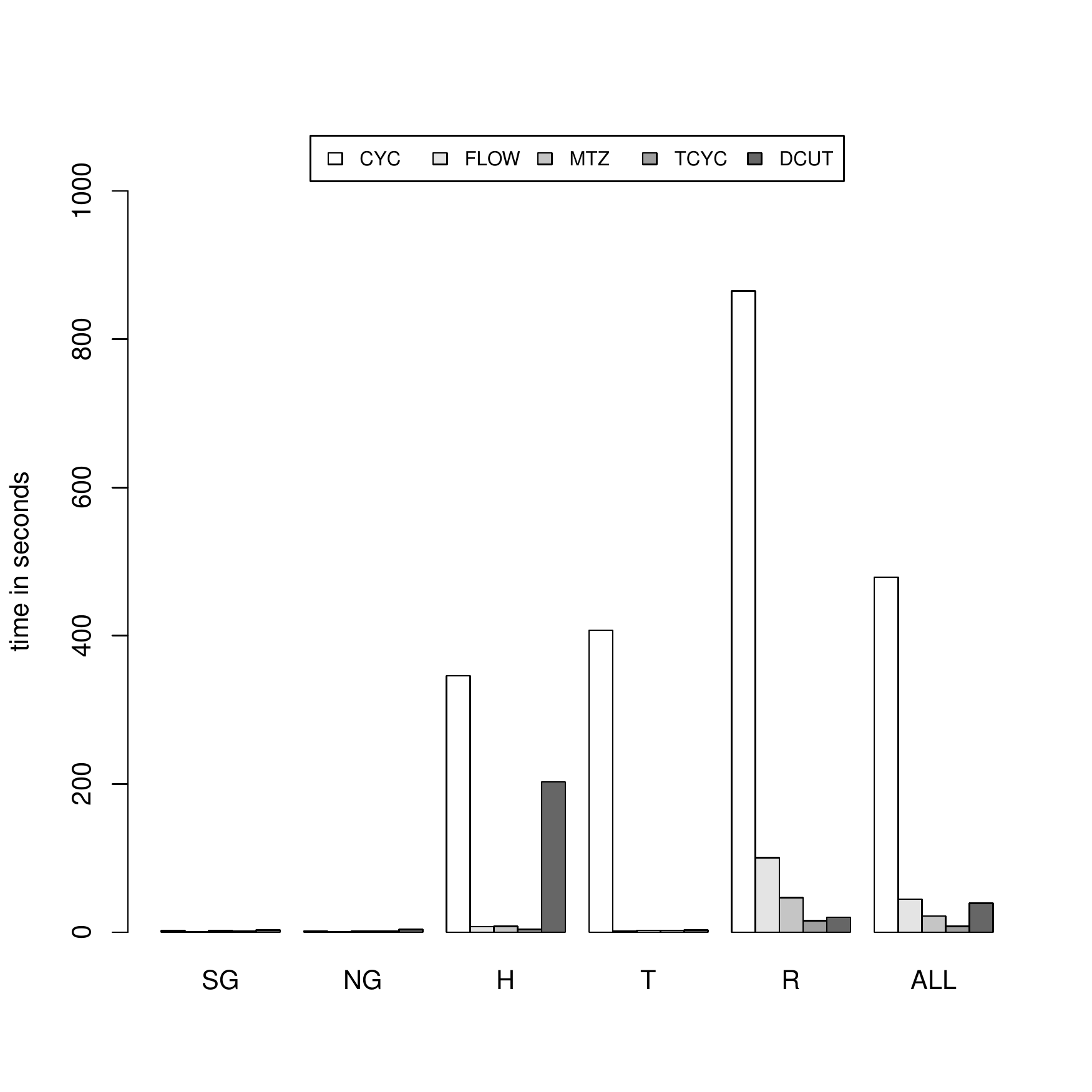}
\end{center}
\caption{Average computation times in seconds over all small instances in each class for each formulation. For each instance, the computation time considered in the averages is either the time to prove optimality (if less than the imposed time limit of 3600 seconds) or this imposed time limit.}
\label{fig:comparison-small-instances}
\end{figure}

The plots in Figure~\ref{fig:comparison-small-instances} compare the average computation times over all individual small instances in each class for each formulation.
The plots show that TCYC is, overall, the best performing formulation for this set of instances as it takes lower average computation times to solve the instances to optimality. Note that TCYC outperforms all other formulations for the hypercube and random graphs. 


\subsection{Large instances}
\label{sec:largeinstances}

Tables~\ref{tab:gridlargeinitsol}-\ref{tab:randlargeinitsol} summarize the results for the large instances.
Tables~\ref{tab:gridlargeinitsol}~and~\ref{tab:gridnqlargeinitsol} show that for the large grid instances, FLOW clearly outperformed the other approaches, achieving 50 and 54 solved instances for the square and non-square grid instances, respectively.
Table~\ref{tab:ipercubelargeinitsol} displays that, for the large hypercube instances, MTZ outperforms the other approaches in terms of number of solutions solved to optimality. This is partly due to the fact that the solver is able to generate good quality feasible solutions earlier in the search process.
Table~\ref{tab:torolargeinitsol} shows that, for the large toroidal instances, FLOW outperformed the other approaches when it comes to the number of solved instances, closely followed by MTZ. 
Table~\ref{tab:randlargeinitsol} displays that, for the large random instances, TCYC presented outstanding results when compared to the other approaches, especially on the average time to solve instances to optimality. Even though MTZ proved optimality of one additional instance, the average times achieved by TCYC and DCUT were much lower. We note, though, that for the larger instances, the open gap is still large, even for the best performing approaches. For very large random instance groups, the open gap was usually very high, with DCUT often achieving the lowest values.

\begin{landscape}

\begin{table}[H]
\centering
\scriptsize
\caption{Results comparing the formulations for the large square grid instances.}\label{tab:gridlargeinitsol}
\begin{tabular}{l rrrr rrrr rrrr rrrr rrrr}
  \hline
  & \multicolumn{4}{c}{CYC} &  \multicolumn{4}{c}{FLOW} & \multicolumn{4}{c}{MTZ} & \multicolumn{4}{c}{TCYC} & \multicolumn{4}{c}{DCUT} \\ \cmidrule[0.4pt](lr{0.125em}){2-5} \cmidrule[0.4pt](lr{0.125em}){6-9} \cmidrule[0.4pt](lr{0.125em}){10-13}
  \cmidrule[0.4pt](lr{0.125em}){14-17}
  \cmidrule[0.4pt](lr{0.125em}){18-21}
Instance & glr & opt & time & gap & glr & opt & time & gap & glr & opt & time & gap & glr & opt & time & gap & glr & opt & time & gap \\
group & (\%)  &  & (s) & (\%) & (\%)  &  & (s) & (\%) & (\%)  &  & (s) & (\%) & (\%)  &  & (s) & (\%) & (\%)  &  & (s) & (\%) \\ \hline
G\_10\_10\_10\_25  & 8.6 & 5 & 71.2 &    & 0.3 & 5 & 3.6 &    & 0.8 & 5 & 33.8 &    & 0.8 & 5 &  \textbf{3.1}  &    & 0.2 & 5 & 5.6 &   \\ 
  G\_10\_10\_10\_50  & 7.1 & 5 & 5.1 &    & 0.5 & 5 & 6.1 &    & 1.1 & 5 & 8.1 &    & 1.0 & 5 &  \textbf{1.4}  &    & 0.3 & 5 & 10.9 &   \\ 
  G\_10\_10\_10\_75  & 8.5 & 5 & 6.2 &    & 0.4 & 5 & 5.4 &    & 0.9 & 5 & 13.5 &    & 1.3 & 5 &  \textbf{1.8}  &    & 0.4 & 5 & 37.5 &   \\ 
  G\_14\_14\_10\_25  & 20.5 & 0 &    & 3.8 & 0.1 & 5 &  \textbf{32.2}  &    & 0.4 & 3 & 1058.0 & 0.1 & 0.4 & 5 & 1009.5 &    & 0.1 & 4 & 95.0 &  0.1 \\ 
  G\_14\_14\_10\_50  & 18.8 & 0 &    & 1.9 & 0.4 & 5 &  \textbf{152.5}  &    & 0.7 & 5 & 408.5 &    & 0.7 & 5 & 551.5 &    & 0.3 & 1 & 2001.1 &  0.3 \\ 
  G\_14\_14\_10\_75  & 19.2 & 0 &    & 2.9 & 0.3 &  \textbf{5}  & 64.6 &    & 0.5 & 3 & 484.9 & 0.1 & 0.8 & 1 & 75.8 & 0.2 & 0.2 & 2 & 94.5 &  0.3 \\ 
  G\_17\_17\_10\_25  & 21.8 & 0 &    & 5.0 & 0.2 &  \textbf{4}  & 237.5 & $<$0.1 & 0.4 & 1 & 2975.9 & 0.1 & 0.4 & 0 &    & 0.2 & 0.2 & 1 & 233.9 &  0.4 \\ 
  G\_17\_17\_10\_50  & 19.8 & 0 &    & 3.3 & 0.4 & 4 & 688.1 & 0.1 & 0.5 &  \textbf{5}  & 1336.3 &    & 0.6 &  \textbf{5}  & 775.9 &    & 0.3 & 1 & 847.3 &  0.3 \\ 
  G\_17\_17\_10\_75  & 20.9 & 0 &    & 4.0 & 0.2 &  \textbf{5}  & 828.0 &    & 0.4 & 1 & 651.8 & 0.1 & 0.6 & 1 & 3575.4 & 0.4 & 0.2 & 0 &    &  0.3 \\ 
  G\_20\_20\_10\_25  & 23.9 & 0 &    & 5.9 & 0.1 &  \textbf{3}  & 1837.5 & $<$0.1 & 0.3 & 0 &    & 0.1 & 0.3 & 0 &    & 0.3 & 0.2 & 0 &    &  0.4 \\ 
  G\_20\_20\_10\_50  & 22.0 & 0 &    & 4.1 & 0.4 &  \textbf{1}  & 2840.8 & 0.1 & 0.5 & 0 &    & 0.1 & 0.8 & 0 &    & 0.3 & 0.3 & 0 &    &  0.5 \\ 
  G\_20\_20\_10\_75  & 22.0 & 0 &    & 4.6 & 0.2 &  \textbf{1}  & 2554.7 & $<$0.1 & 0.3 & 0 &    & 0.1 & 0.6 & 0 &    & 0.5 & 0.2 & 0 &    &  0.4 \\ 
  G\_23\_23\_10\_25  & 24.8 & 0 &    & 6.1 & 0.1 &  \textbf{1}  & 473.5 & 0.1 & 0.2 & 0 &    & 0.2 & 0.2 & 0 &    & 0.3 & 0.1 & 0 &    &  0.3 \\ 
  G\_23\_23\_10\_50  & 23.0 & 0 &    & 4.4 & 0.4 & 0 &    &  \textbf{0.1}  & 0.4 & 0 &    & 0.2 & 0.7 & 0 &    & 0.4 & 0.3 & 0 &    &  0.6 \\ 
  G\_23\_23\_10\_75  & 23.8 & 0 &    & 5.0 & 0.2 & 0 &    &  \textbf{0.1}  & 0.3 & 0 &    &  \textbf{0.1}  & 0.5 & 0 &    & 0.5 & 0.2 & 0 &    &  0.5 \\ \hline
  Average  & 19.0 &    & 27.5 & 4.2 & 0.3 &    & 748.1 &  {0.1}  & 0.5 &    & 774.5 &  {0.1}  & 0.7 &    & 749.3 & 0.3 & 0.2 &    & 415.7 &  0.4 \\
  Total &  & 15 &  &  &  & 49 &  &  &  & 33 &  &  &  & 32 &  &  &  & 24 &  &  \\ 
   \hline
\end{tabular}
\end{table}

\begin{table}[H]
\centering
\scriptsize
\caption{Results comparing the formulations for the large non-square grid instances.}\label{tab:gridnqlargeinitsol}
\begin{tabular}{l rrrr rrrr rrrr rrrr rrrr}
  \hline
  & \multicolumn{4}{c}{CYC} &  \multicolumn{4}{c}{FLOW} & \multicolumn{4}{c}{MTZ} & \multicolumn{4}{c}{TCYC} & \multicolumn{4}{c}{DCUT} \\ \cmidrule[0.4pt](lr{0.125em}){2-5} \cmidrule[0.4pt](lr{0.125em}){6-9} \cmidrule[0.4pt](lr{0.125em}){10-13}
  \cmidrule[0.4pt](lr{0.125em}){14-17}
  \cmidrule[0.4pt](lr{0.125em}){18-21}
Instance & glr & opt & time & gap & glr & opt & time & gap & glr & opt & time & gap & glr & opt & time & gap & glr & opt & time & gap \\
group & (\%)  &  & (s) & (\%) & (\%)  &  & (s) & (\%) & (\%)  &  & (s) & (\%) & (\%)  &  & (s) & (\%) & (\%)  &  & (s) & (\%) \\ \hline
GNQ\_13\_7\_10\_25  & 7.3 & 5 & 15.6 &    & 0.4 & 5 & 4.7 &    & 0.9 & 5 & 23.4 &    & 0.9 & 5 &  \textbf{1.8}  &    & 0.3 & 5 & 250.7 &   \\ 
  GNQ\_13\_7\_10\_50  & 5.1 & 5 & 2.3 &    & 0.4 & 5 &  \textbf{1.2}  &    & 1.0 & 5 & 1.9 &    & 0.8 & 5 & 1.3 &    & 0.2 & 5 & 8.1 &   \\ 
  GNQ\_13\_7\_10\_75  & 5.7 & 5 & 3.6 &    & 0.4 & 5 & 2.7 &    & 1.0 & 5 & 11.3 &    & 1.0 & 5 &  \textbf{1.5}  &    & 0.2 & 5 & 7.1 &   \\ 
  GNQ\_18\_11\_10\_25  & 18.7 & 0 &    & 3.4 & 0.3 &  \textbf{5}  & 202.3 &    & 0.5 & 4 & 1253.5 & 0.2 & 0.5 & 4 & 227.2 & 0.2 & 0.2 & 1 & 77.7 &  0.2 \\ 
  GNQ\_18\_11\_10\_50  & 16.5 & 0 &    & 2.2 & 0.4 & 5 & 162.0 &    & 0.6 & 5 & 411.4 &    & 0.7 & 5 &  \textbf{150.2}  &    & 0.2 & 3 & 120.2 &  0.4 \\ 
  GNQ\_18\_11\_10\_75  & 17.0 & 0 &    & 2.6 & 0.3 &  \textbf{5}  & 148.7 &    & 0.5 & 3 & 1237.7 & 0.1 & 0.7 & 3 & 544.0 & 0.4 & 0.2 & 2 & 939.8 &  0.2 \\ 
  GNQ\_23\_13\_10\_25  & 20.7 & 0 &    & 5.4 & 0.1 &  \textbf{5}  & 738.5 &    & 0.3 & 0 &    & 0.1 & 0.3 & 0 &    & 0.2 & 0.2 & 1 & 1515.2 &  0.2 \\ 
  GNQ\_23\_13\_10\_50  & 18.7 & 0 &    & 3.5 & 0.4 &  \textbf{4}  & 905.8 & 0.1 & 0.6 & 4 & 1818.8 & 0.1 & 0.8 & 2 & 1358.1 & 0.3 & 0.3 & 0 &    &  0.5 \\ 
  GNQ\_23\_13\_10\_75  & 19.2 & 0 &    & 4.1 & 0.2 &  \textbf{4}  & 1194.7 & 0.1 & 0.4 & 1 & 1579.1 & 0.1 & 0.6 & 0 &    & 0.3 & 0.2 & 0 &    &  0.3 \\ 
  GNQ\_26\_15\_10\_25  & 22.4 & 0 &    & 5.8 & 0.1 &  \textbf{5}  & 1690.5 &    & 0.2 & 0 &    & 0.1 & 0.3 & 0 &    & 0.2 & 0.1 & 0 &    &  0.3 \\ 
  GNQ\_26\_15\_10\_50  & 20.0 & 0 &    & 3.6 & 0.4 &  \textbf{2}  & 1995.1 & 0.1 & 0.5 & 1 & 1533.7 & 0.1 & 0.6 & 1 & 532.8 & 0.3 & 0.3 & 1 & 1603.0 &  0.4 \\ 
  GNQ\_26\_15\_10\_75  & 20.8 & 0 &    & 4.6 & 0.2 &  \textbf{4}  & 973.2 & $<$0.1 & 0.3 & 0 &    & 0.1 & 0.5 & 0 &    & 0.4 & 0.1 & 0 &    &  0.3 \\ 
  GNQ\_29\_17\_10\_25  & 23.4 & 0 &    & 6.2 & 0.1 &  \textbf{3}  & 2371.0 & 0.1 & 0.2 & 0 &    & 0.1 & 0.2 & 0 &    & 0.2 & 0.1 & 0 &    &  0.2 \\ 
  GNQ\_29\_17\_10\_50  & 21.6 & 0 &    & 4.4 & 0.4 & 0 &    &  \textbf{0.1}  & 0.6 & 0 &    & 0.2 & 0.7 & 0 &    & 0.4 & 0.3 & 0 &    &  0.7 \\ 
  GNQ\_29\_17\_10\_75  & 22.3 & 0 &    & 4.7 & 0.2 & 0 &    &  \textbf{0.1}  & 0.3 & 0 &    & 0.2 & 0.5 & 0 &    & 0.4 & 0.2 & 0 &    &  0.4 \\ \hline
  Average  & 17.3 &    & 7.2 & 4.2 & 0.3 &    & 799.3 &  \textbf{0.1}  & 0.5 &    & 874.5 &  \textbf{0.1}  & 0.6 &    & 352.1 & 0.3 & 0.2 &    & 565.2 &  0.3 \\
  Total &  & 15 &  &  &  & 57 &  &  &  & 33 &  &  &  & 30 &  &  &  & 23 &  &  \\ 
   \hline
\end{tabular}
\end{table}

\begin{table}[H]
\centering
\scriptsize
\caption{Results comparing the formulations for the large hypercube instances.}\label{tab:ipercubelargeinitsol}
\begin{tabular}{l rrrr rrrr rrrr rrrr rrrr}
  \hline
  & \multicolumn{4}{c}{CYC} &  \multicolumn{4}{c}{FLOW} & \multicolumn{4}{c}{MTZ} & \multicolumn{4}{c}{TCYC} & \multicolumn{4}{c}{DCUT} \\ \cmidrule[0.4pt](lr{0.125em}){2-5} \cmidrule[0.4pt](lr{0.125em}){6-9} \cmidrule[0.4pt](lr{0.125em}){10-13}
  \cmidrule[0.4pt](lr{0.125em}){14-17}
  \cmidrule[0.4pt](lr{0.125em}){18-21}
Instance & glr & opt & time & gap & glr & opt & time & gap & glr & opt & time & gap & glr & opt & time & gap & glr & opt & time & gap \\
group & (\%)  &  & (s) & (\%) & (\%)  &  & (s) & (\%) & (\%)  &  & (s) & (\%) & (\%)  &  & (s) & (\%) & (\%)  &  & (s) & (\%) \\ \hline
H\_7\_10\_25 &29.5&0&  &13.2&5.0&0&  &1.9&5.1& \textbf{1} &2489.6&1.4&5.9&0&  &1.4&4.8&0&  & 2.8 \\ 
  H\_7\_10\_50 &27.1&0&  &5.8&2.8&5&317.4&  &2.8&5& \textbf{135.8} &  &4.4&5&222.8&  &2.6&4&63.8& 1.1 \\ 
  H\_7\_10\_75 &23.5&0&  &3.8&1.7&5&16.6&  &1.8&5& \textbf{7.1} &  &3.3&5&26.9&  &1.5&5&42.1&  \\ 
  H\_8\_10\_25 &31.1&0&  &18.5&3.8&0&  &3.4&3.8&0&  & \textbf{2.3} &4.7&0&  &3.6&3.7&0&  & 3.1 \\ 
  H\_8\_10\_50 &29.6&0&  &13.6&3.6&0&  &2.4&3.6&0&  & \textbf{1.6} &5.4&0&  &3.1&3.6&0&  & 2.2 \\ 
  H\_8\_10\_75 &26.8&0&  &10.2&1.9&1&654.4&0.6&1.9& \textbf{4} &451.4&0.5&3.7&0&  &1.5&1.8& \textbf{4} &1398.5& 1.1 \\ 
  H\_9\_10\_25 &34.3&0&  &23.1&5.5&0&  &5.7&5.6&0&  & \textbf{5.3} &6.6&0&  &5.7&5.5&0&  & 5.6 \\ 
  H\_9\_10\_50 &31.9&0&  &17.8&4.3&0&  &4.3&4.4&0&  & \textbf{3.9} &6.0&0&  &5&4.3&0&  & 4.0 \\ 
  H\_9\_10\_75 &30.0&0&  &15.3&3.2&0&  &2.8&3.1&0&  & \textbf{2.4} &5.1&0&  &3.9&3.1&0&  & \textbf{2.4} \\ \hline
  Average &29.3&  &  &13.5&3.6&  &329.5&3&3.6&  &771& {2.5} &5.0&  &124.9&3.5&3.4&  &501.5& 2.8 \\ 
  Total &  & 0 &  &  &  & 11 &  &  &  & 15 &  &  &  & 10 &  &  &  & 13 &  &  \\ 
   \hline
\end{tabular}
\end{table}

\begin{table}[H]
\centering
\scriptsize
\caption{Results comparing the formulations for the large toroidal instances.}\label{tab:torolargeinitsol}
\begin{tabular}{l rrrr rrrr rrrr rrrr rrrr}
  \hline
  & \multicolumn{4}{c}{CYC} &  \multicolumn{4}{c}{FLOW} & \multicolumn{4}{c}{MTZ} & \multicolumn{4}{c}{TCYC} & \multicolumn{4}{c}{DCUT} \\ \cmidrule[0.4pt](lr{0.125em}){2-5} \cmidrule[0.4pt](lr{0.125em}){6-9} \cmidrule[0.4pt](lr{0.125em}){10-13}
  \cmidrule[0.4pt](lr{0.125em}){14-17}
  \cmidrule[0.4pt](lr{0.125em}){18-21}
Instance & glr & opt & time & gap & glr & opt & time & gap & glr & opt & time & gap & glr & opt & time & gap & glr & opt & time & gap \\
group & (\%)  &  & (s) & (\%) & (\%)  &  & (s) & (\%) & (\%)  &  & (s) & (\%) & (\%)  &  & (s) & (\%) & (\%)  &  & (s) & (\%) \\ \hline
T\_10\_10\_10\_25  & 15.2 & 0 &    & 2.9 & 0.8 & 5 & 8.5 &    & 1.1 & 5 & 34.8 &    & 1.2 & 5 & 19.4 &    & 0.6 & 5 &  \textbf{7.9}  &   \\ 
  T\_10\_10\_10\_50  & 10.7 & 5 & 1040.0 &    & 0.7 & 5 &  \textbf{3.4}  &    & 1.2 & 5 & 15.3 &    & 1.5 & 5 & 4.2 &    & 0.5 & 5 & 100.4 &   \\ 
  T\_10\_10\_10\_75  & 8.7 & 5 & 69.5 &    & 0.4 & 5 &  \textbf{1.4}  &    & 1.0 & 5 & 1.8 &    & 1.1 & 5 & 2.1 &    & 0.5 & 4 & 3.9 &  0.6 \\ 
  T\_14\_14\_10\_25  & 21.4 & 0 &    & 4.9 & 0.5 & 5 &  \textbf{407.9}  &    & 0.7 & 5 & 676.2 &    & 0.9 & 1 & 11.1 & 0.5 & 0.5 & 5 & 453.7 &   \\ 
  T\_14\_14\_10\_50  & 16.6 & 0 &    & 2.1 & 0.5 &  \textbf{5}  & 50.5 &    & 0.7 &  \textbf{5}  & 94.5 &    & 1.1 & 3 & 110.6 & 0.3 & 0.3 & 4 & 555.9 &  0.3 \\ 
  T\_14\_14\_10\_75  & 13.0 & 1 & 1194.3 & 1.1 & 0.5 & 5 &  \textbf{53.4}  &    & 0.9 & 5 & 55.2 &    & 0.7 & 5 & 64.0 &    & 0.2 & 0 &    &  0.5 \\ 
  T\_17\_17\_10\_25  & 20.5 & 0 &    & 5.5 & 0.7 &  \textbf{2}  & 942.0 & 0.1 & 0.8 & 1 & 287.9 & 0.1 & 1.0 & 0 &    & 0.6 & 0.5 &  \textbf{2}  & 541.6 &  0.5 \\ 
  T\_17\_17\_10\_50  & 18.6 & 0 &    & 2.9 & 0.6 &  \textbf{5}  & 500.2 &    & 0.7 &  \textbf{5}  & 762.7 &    & 0.9 & 0 &    & 0.4 & 0.3 & 1 & 1285.2 &  0.6 \\ 
  T\_17\_17\_10\_75  & 16.3 & 0 &    & 1.8 & 0.6 &  \textbf{5}  & 723.2 &    & 0.8 &  \textbf{5}  & 403.8 &    & 0.8 & 2 & 1424.9 & 0.3 & 0.3 & 0 &    &  0.8 \\
  T\_20\_20\_10\_25  & 22.8 & 0 &    & 6.2 & 0.6 &  \textbf{1}  & 3461.2 & 0.2 & 0.7 & 0 &    & 0.2 & 0.9 & 0 &    & 0.8 & 0.4 & 0 &    &  0.6 \\ 
  T\_20\_20\_10\_50  & 19.1 & 0 &    & 3.2 & 0.5 &  \textbf{4}  & 2302.7 & 0.1 & 0.6 & 2 & 2144.7 & 0.1 & 0.9 & 0 &    & 0.5 & 0.3 & 0 &    &  0.6 \\ 
  T\_20\_20\_10\_75  & 15.6 & 0 &    & 2.0 & 0.5 &  \textbf{5}  & 1447.6 &    & 0.7 & 4 & 886.3 & 0.1 & 0.8 & 0 &    & 0.3 & 0.2 & 0 &    &  0.6 \\ 
  T\_23\_23\_10\_25  & 24.7 & 0 &    & 6.5 & 0.6 & 0 &    & 0.3 & 0.7 & 0 &    & 0.4 & 0.9 & 0 &    & 0.9 & 0.5 & 0 &    &  0.6 \\ 
  T\_23\_23\_10\_50  & 18.5 & 0 &    & 3.2 & 0.5 &  \textbf{1}  & 2526.1 & 0.1 & 0.6 &  \textbf{1}  & 3105.0 & 0.1 & 0.8 & 0 &    & 0.5 & 0.3 & 0 &    &  0.8 \\ 
  T\_23\_23\_10\_75  & 17.0 & 0 &    & 2.2 & 0.4 & 1 & 2847.3 & 0.1 & 0.7 &  \textbf{2}  & 2359.1 & 0.1 & 0.8 & 0 &    & 0.4 & 0.3 & 0 &    &  0.8 \\ \hline
  Average  & 17.2 &    & 767.9 & 3.4 & 0.6 &    & 1091.1 & 0.1 & 0.8 &    & 832.9 & 0.2 & 1.0 &    & 233.8 & 0.5 & 0.4 &    & 421.2 &  0.6 \\
  Total &  & 11 &  &  &  & 54 &  &  &  & 50 &  &  &  & 26 &  &  &  & 26 &  &  \\ 
   \hline
\end{tabular}
\end{table}

\begin{table}[H]
\centering
\scriptsize
\caption{Results comparing the formulations for the large random instances.}\label{tab:randlargeinitsol}
\begin{tabular}{l rrrr rrrr rrrr rrrr rrrr}
  \hline
  & \multicolumn{4}{c}{CYC} &  \multicolumn{4}{c}{FLOW} & \multicolumn{4}{c}{MTZ} & \multicolumn{4}{c}{TCYC} & \multicolumn{4}{c}{DCUT} \\ \cmidrule[0.4pt](lr{0.125em}){2-5} \cmidrule[0.4pt](lr{0.125em}){6-9} \cmidrule[0.4pt](lr{0.125em}){10-13}
  \cmidrule[0.4pt](lr{0.125em}){14-17}
  \cmidrule[0.4pt](lr{0.125em}){18-21}
Instance & glr & opt & time & gap & glr & opt & time & gap & glr & opt & time & gap & glr & opt & time & gap & glr & opt & time & gap \\
group & (\%)  &  & (s) & (\%) & (\%)  &  & (s) & (\%) & (\%)  &  & (s) & (\%) & (\%)  &  & (s) & (\%) & (\%)  &  & (s) & (\%) \\ \hline  &  & (s) & (\%) &  &  & (s) & (\%) &  &  & (s) & (\%) \\
  \hline
  R\_100\_247\_10\_25  & 22.9 & 0 &    & 6.7 & 1.3 & 5 & 6.5 &    & 1.8 & 5 &  \textbf{5.1}  &    & 1.7 & 5 & 11.7 &    & 1.3 & 5 & 5.4 &   \\ 
  R\_100\_247\_10\_50  & 24.3 & 0 &    & 6.0 & 1.7 & 5 &  \textbf{7.0}  &    & 2.1 & 5 &  \textbf{7.0}  &    & 2.4 & 5 & 56.0 &    & 1.6 & 5 & 8.9 &   \\ 
  R\_100\_247\_10\_75  & 24.6 & 0 &    & 5.5 & 1.7 & 5 & 6.4 &    & 2.0 & 5 &  \textbf{4.5}  &    & 2.3 & 5 & 52.3 &    & 1.7 & 5 & 8.3 &   \\ 
  R\_100\_841\_10\_25  & 54.0 & 0 &    & 28.3 & 14.4 & 5 & 1928.0 &    & 14.7 & 5 & 788.4 &    & 11.9 & 5 &  \textbf{117.2}  &  & 14.8 & 5& 1008.9 &   \\  
  R\_100\_841\_10\_50  & 56.7 & 0 &    & 29.8 & 17.3 & 5 & 2521.2 &    & 17.5 & 5 & 961.7 &    & 14.6 & 5 &  \textbf{210.2}  & & 17.7 & 5 & 1022.5 &   \\   
  R\_100\_841\_10\_75  & 55.5 & 0 &    & 28.9 & 17.5 & 5 & 1877.5 &    & 17.6 & 5 & 759.5 &    & 14.6 & 5 &  \textbf{181.2}  &  & 17.8 & 5 & 826.8 & \\
  R\_100\_3069\_10\_25  & 81.9 & 0 &    & 66.4 & 64.8 & 0 &    & 34.6 & 64.8 & 5 & 2589.2 &    & 62.6 & 5 &  \textbf{827.4}  &  & 41.5 & 5& 1004.4 &   \\ 
  R\_100\_3069\_10\_50  & 81.7 & 0 &    & 64.7 & 63.8 & 1 & 3190.4 & 27.7 & 64.4 & 5 & 1879.3 &    & 61.2 & 5 & 683.3 &  & 39.2 & 5 &  \textbf{557.4}  &   \\  
  R\_100\_3069\_10\_75  & 81.9 & 0 &    & 63.7 & 63.8 & 1 & 2817.0 & 23.6 & 64.0 & 5 & 1738.9 &    & 58.0 & 5 & 602.9 & & 38.7 & 5 & \textbf{453.6}  &   \\ 
  R\_200\_796\_10\_25  & 39.7 & 0 &    & 23.9 & 6.0 & 0 &    & 4.2 & 6.4 & 0 &    &  \textbf{2.8}  & 5.7 & 0 &    & 2.9 & 6.1 & 0 &    &  3.2 \\ 
  R\_200\_796\_10\_50  & 39.4 & 0 &    & 22.8 & 6.3 & 0 &    & 4.0 & 6.3 & 0 &    &  \textbf{2.2}  & 5.8 & 0 &    & 2.9 & 6.1 & 0 &    &  2.7 \\ 
  R\_200\_796\_10\_75  & 38.8 & 0 &    & 21.4 & 6.3 & 0 &    & 3.3 & 6.4 & 1 & 2751.6 &  \textbf{1.7}  & 5.6 & 0 &    & 2.7 & 6.6 & 0 &    &  2.4 \\ 
  R\_200\_3184\_10\_25  & 69.0 & 0 &    & 55.9 & 34.5 & 0 &    & 33.2 & 34.5 & 0 &    & 27.5 & 32.0 & 0 &    &  \textbf{26.7}  & 34.8 & 0 &    &  30.0 \\ 
  R\_200\_3184\_10\_50  & 69.4 & 0 &    & 56.4 & 35.0 & 0 &    & 34.1 & 35.2 & 0 &    &  \textbf{27.3}  & 33.5 & 0 &    & 27.6 & 35.3 & 0 &    &  30.4 \\ 
  R\_200\_3184\_10\_75  & 69.2 & 0 &    & 55.4 & 34.2 & 0 &    & 32.8 & 34.4 & 0 &    & 26.2 & 32.5 & 0 &    &  \textbf{25.7}  & 34.8 & 0 &    &  28.9 \\ 
  R\_200\_12139\_10\_25  & 89.0 & 0 &    & 83.3 & 84.0 & 0 &    & 82.9 & 78.8 & 0 &    & 78.3 & 75.7 & 0 &    & 74.1 & 58.2 & 0 &    &  \textbf{54.7} \\ 
  R\_200\_12139\_10\_50  & 88.4 & 0 &    & 81.9 & 80.7 & 0 &    & 80.8 & 76.4 & 0 &    & 75.8 & 74.0 & 0 &    & 72.4 & 55.2 & 0 &    &  \textbf{50.9} \\ 
  R\_200\_12139\_10\_75  & 87.8 & 0 &    & 80.9 & 82.0 & 0 &    & 81.5 & 75.6 & 0 &    & 74.6 & 74.2 & 0 &    & 70.6 & 52.7 & 0 &    &  \textbf{48.6} \\ 
  R\_300\_1644\_10\_25  & 45.6 & 0 &    & 30.4 & 9.0 & 0 &    & 8.5 & 9.2 & 0 &    & 7.5 & 8.2 & 0 &    &  \textbf{6.4}  & 9.3 & 0 &    &  7.3 \\ 
  R\_300\_1644\_10\_50  & 48.6 & 0 &    & 34.1 & 11.6 & 0 &    & 11.2 & 11.8 & 0 &    & 10.1 & 11.2 & 0 &    &  \textbf{9.1}  & 12.2 & 0 &    &  10.0 \\ 
  R\_300\_1644\_10\_75  & 47.6 & 0 &    & 33.5 & 11.4 & 0 &    & 10.8 & 11.6 & 0 &    & 9.8 & 10.9 & 0 &    &  \textbf{8.8}  & 11.9 & 0 &    &  10.0 \\ 
  R\_300\_7026\_10\_25  & 76.3 & 0 &    & 67.8 & 48.8 & 0 &    & 52.4 & 47.3 & 0 &    & 46.7 & 45.7 & 0 &    &  \textbf{43.1}  & 47.9 & 0 &    &  46.6 \\ 
  R\_300\_7026\_10\_50  & 76.5 & 0 &    & 67.8 & 48.8 & 0 &    & 51.6 & 47.6 & 0 &    & 46.8 & 45.8 & 0 &    &  \textbf{42.2}  & 48.2 & 0 &    &  46.9 \\ 
  R\_300\_7026\_10\_75  & 76.4 & 0 &    & 67.2 & 48.1 & 0 &    & 50.8 & 47.2 & 0 &    & 46.4 & 45.5 & 0 &    &  \textbf{41.9}  & 48.0 & 0 &    &  46.7 \\ 
  R\_300\_27209\_10\_25  & 91.9 & 0 &    & 88.8 & 94.2 & 0 &    & 94.2 & 88.5 & 0 &    & 88.5 & 81.7 & 0 &    & 81.6 & 66.9 & 0 &    &  \textbf{67.4} \\ 
  R\_300\_27209\_10\_50  & 90.7 & 0 &    & 86.7 & 90.7 & 0 &    & 90.7 & 86.3 & 0 &    & 86.3 & 78.8 & 0 &    & 77.4 & 61.3 & 0 &    &  \textbf{61.4} \\ 
  R\_300\_27209\_10\_75  & 91.4 & 0 &    & 88.1 & 92.7 & 0 &    & 92.7 & 87.6 & 0 &    & 87.6 & 80.6 & 0 &    & 79.4 & 64.3 & 0 &    &  \textbf{64.2} \\ 
  R\_400\_2793\_10\_25  & 51.8 & 0 &    & 39.8 & 14.3 & 0 &    & 14.4 & 14.5 & 0 &    & 14.2 & 13.6 & 0 &    &  \textbf{12.0}  & 14.4 & 0 &    &  13.8 \\ 
  R\_400\_2793\_10\_50  & 55.2 & 0 &    & 42.6 & 17.6 & 0 &    & 17.6 & 17.5 & 0 &    & 17.3 & 16.9 & 0 &    &  \textbf{15.3}  & 17.9 & 0 &    &  17.0 \\ 
  R\_400\_2793\_10\_75  & 54.8 & 0 &    & 42.2 & 16.7 & 0 &    & 16.7 & 16.9 & 0 &    & 16.5 & 16.3 & 0 &    &  \textbf{14.7}  & 17.0 & 0 &    &  16.3 \\ 
  R\_400\_12369\_10\_25  & 80.3 & 0 &    & 76.2 & 68.0 & 0 &    & 68.1 & 56.7 & 0 &    & 59.4 & 54.1 & 0 &    &  \textbf{52.8}  & 56.5 & 0 &    &  56.4 \\ 
  R\_400\_12369\_10\_50  & 78.0 & 0 &    & 72.3 & 63.8 & 0 &    & 64.1 & 50.8 & 0 &    & 54.1 & 49.0 & 0 &    &  \textbf{47.1}  & 52.3 & 0 &    &  51.3 \\ 
  R\_400\_12369\_10\_75  & 79.4 & 0 &    & 74.4 & 66.2 & 0 &    & 66.5 & 54.0 & 0 &    & 60.0 & 52.1 & 0 &    &  \textbf{50.5}  & 54.2 & 0 &    &  53.9 \\ 
  R\_400\_48279\_10\_25  & 93.6 & 0 &    & 92.3 & 95.4 & 0 &    & 95.4 & 95.4 & 0 &    & 95.4 & 86.2 & 0 &    &  \textbf{85.0}  & 90.8 & 0 &    &  90.8 \\ 
  R\_400\_48279\_10\_50  & 92.1 & 0 &    & 90.4 & 94.1 & 0 &    & 94.1 & 94.1 & 0 &    & 94.1 & 81.7 & 0 &    &  \textbf{80.1}  & 86.6 & 0 &    &  88.3 \\ 
  R\_400\_48279\_10\_75  & 93.3 & 0 &    & 92.1 & 95.2 & 0 &    & 95.2 & 95.2 & 0 &    & 95.2 & 81.7 & 0 &    &  \textbf{83.3}  & 86.8 & 0 &    &  90.4 \\ 
  R\_500\_4241\_10\_25  & 57.3 & 0 &    & 49.8 & 19.6 & 0 &    & 21.4 & 19.5 & 0 &    & 19.6 & 18.6 & 0 &    &  \textbf{17.3}  & 19.7 & 0 &    &  19.3 \\ 
  R\_500\_4241\_10\_50  & 58.5 & 0 &    & 50.0 & 21.5 & 0 &    & 23.2 & 21.3 & 0 &    & 21.2 & 20.4 & 0 &    &  \textbf{19.3}  & 21.5 & 0 &    &  21.1 \\ 
  R\_500\_4241\_10\_75  & 59.2 & 0 &    & 50.1 & 21.6 & 0 &    & 23.7 & 21.4 & 0 &    & 21.3 & 20.5 & 0 &    &  \textbf{19.3}  & 21.7 & 0 &    &  21.2 \\ 
  R\_500\_19211\_10\_25  & 83.2 & 0 &    & 82.4 & 82.8 & 0 &    & 86.2 & 72.5 & 0 &    & 72.7 & 60.3 & 0 &    &  \textbf{60.3}  & 69.2 & 0 &    &  69.6 \\ 
  R\_500\_19211\_10\_50  & 80.1 & 0 &    & 78.0 & 67.5 & 0 &    & 70.7 & 65.4 & 0 &    & 66.2 & 53.4 & 0 &    &  \textbf{53.2}  & 61.4 & 0 &    &  62.8 \\ 
  R\_500\_19211\_10\_75  & 82.8 & 0 &    & 81.7 & 74.5 & 0 &    & 82.7 & 71.6 & 0 &    & 71.6 & 59.2 & 0 &    &  \textbf{58.9}  & 67.8 & 0 &    &  68.4 \\ 
  R\_500\_75349\_10\_25  & 94.6 &  n.a.  &  n.a.  &  n.a.  & 96.1 & 0 &    & 96.1 & 96.1 & 0 &    & 96.1 & 88.7 & 0 &    &  \textbf{89.4}  & 96.1 & 0 &    &  96.1 \\ 
  R\_500\_75349\_10\_50  & 93.3 &  n.a.  &  n.a.  &  n.a.   & 95.1 & 0 &    & 95.1 & 95.1 & 0 &    & 95.1 & 85.0 & 0 &    &  \textbf{85.7}  & 95.1 & 0 &    &  95.1 \\ 
  R\_500\_75349\_10\_75  & 93.9 &  n.a.  &  n.a.  &  n.a.   & 95.5 & 0 &    & 95.5 & 95.5 & 0 &    & 95.5 & 87.2 & 0 &    &  \textbf{88.7}  & 95.5 & 0 &    &  95.5 \\ \hline
  Average  & 68.9 &    &    & 56.2 & 48.4 &    & 1544.3 & 52.1 & 46.5 &    & 1148.5 & 50.3 & 42.7 &    & 304.7 & 45.2 & 41.3 &    & 544.0 &  45.5 \\
  Total &  & 0 &  &  &  & 32 &  &  &  & 46 &  &  &  & 45 &  &  &  & 45 &  &  \\ 
 
   \hline
\end{tabular}
\end{table}

\end{landscape}

\begin{figure}[H]
\begin{center}
     \includegraphics[width=0.8\textwidth]{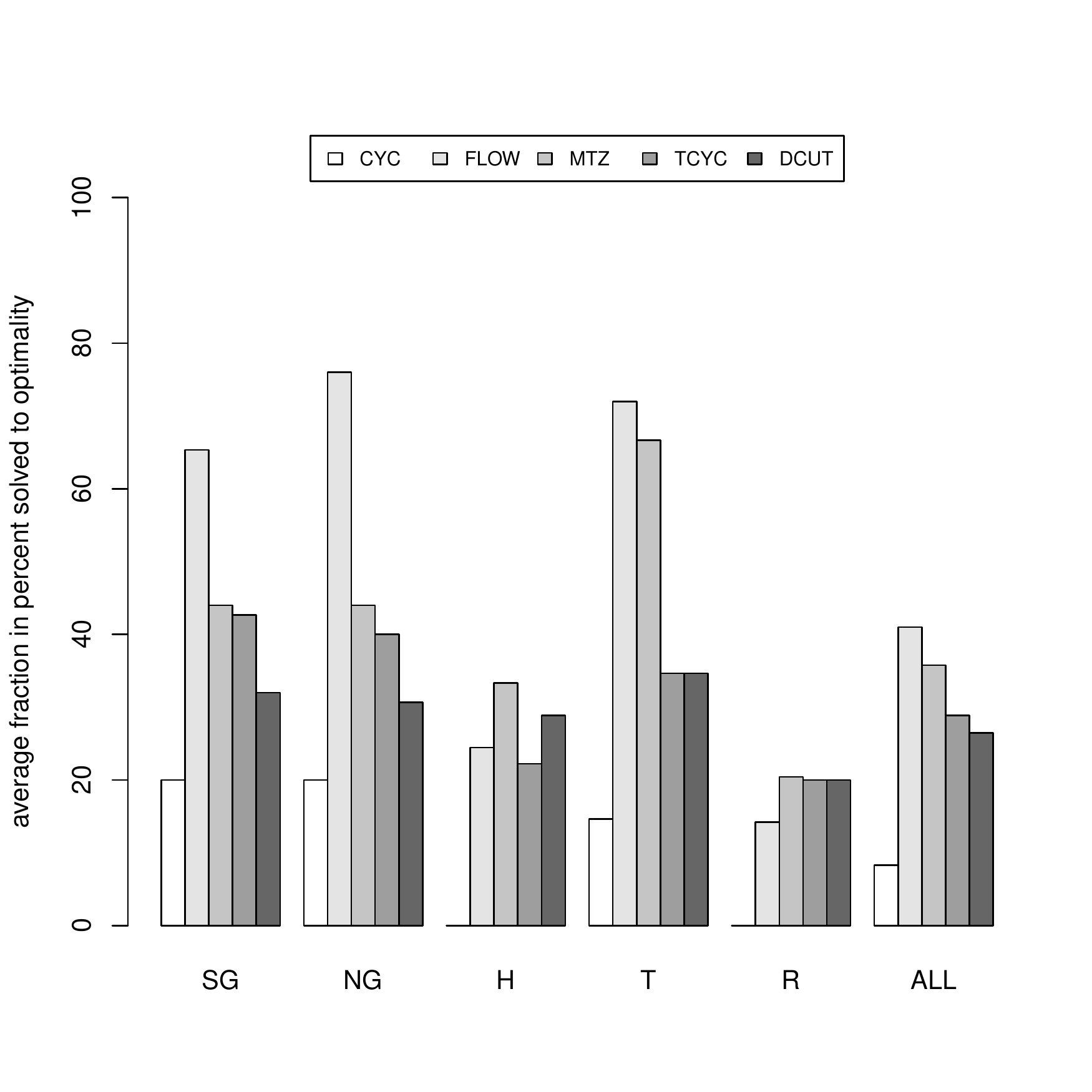}
\end{center}
\caption{Average fraction in percent over all individual large instances solved to optimality in each class for each formulation. For each instance, a time limit of 3600 seconds was imposed.}\label{fig:comparison-large-instances}
\end{figure}

The plots in Figure~\ref{fig:comparison-large-instances}
compare the average fraction in percent over all individual large instances solved to optimality within the time limit of 3600 seconds in each class for each formulation.
The plots show that FLOW performs the best for 
three out of the five classes (square grid, non-square grid, and toroidal).

\begin{figure}[H]
\subfigure[SG]{
   \includegraphics[scale =0.45] {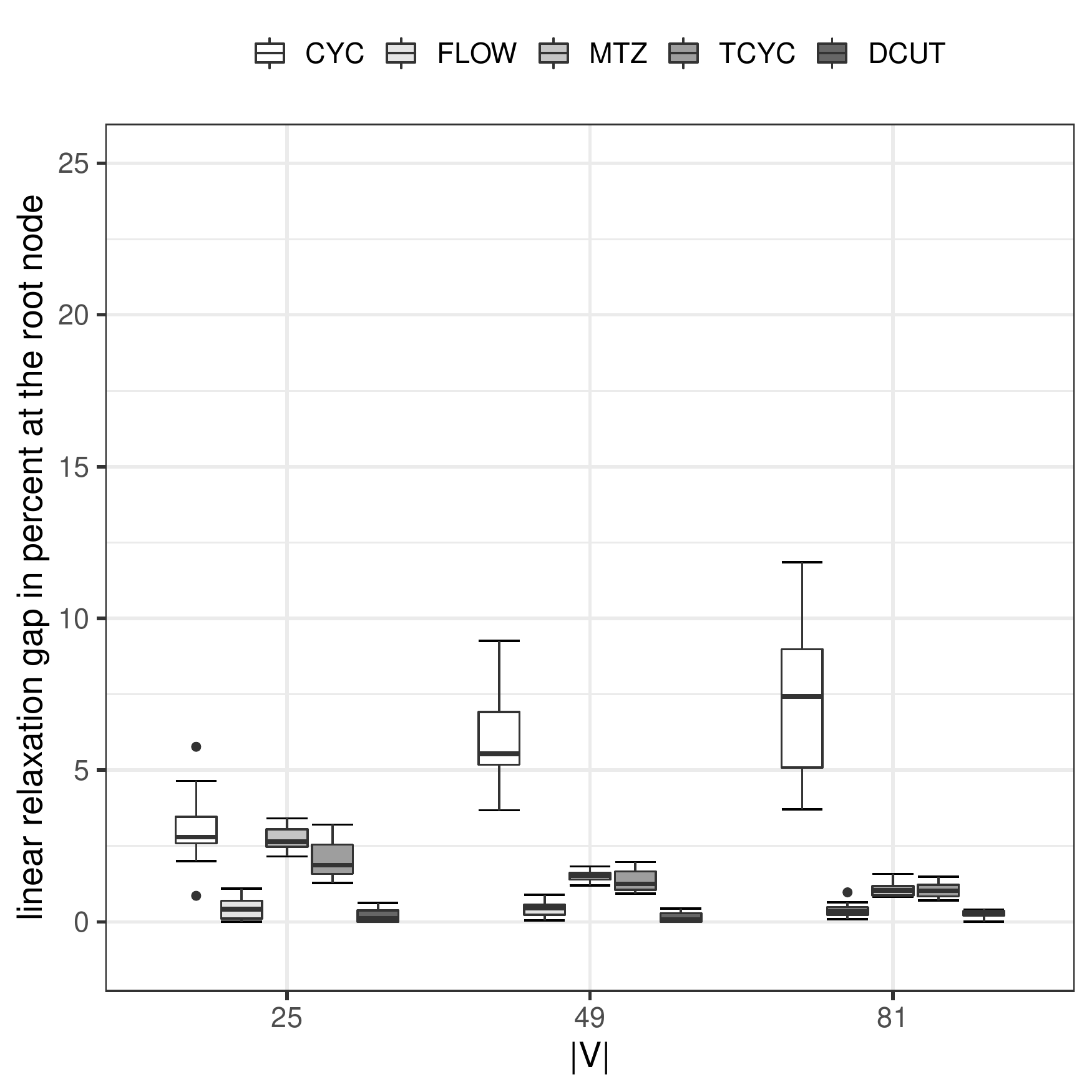}
   \label{fig:glrboxplota}
 }
  \subfigure[NG]{
   \includegraphics[scale =0.45] {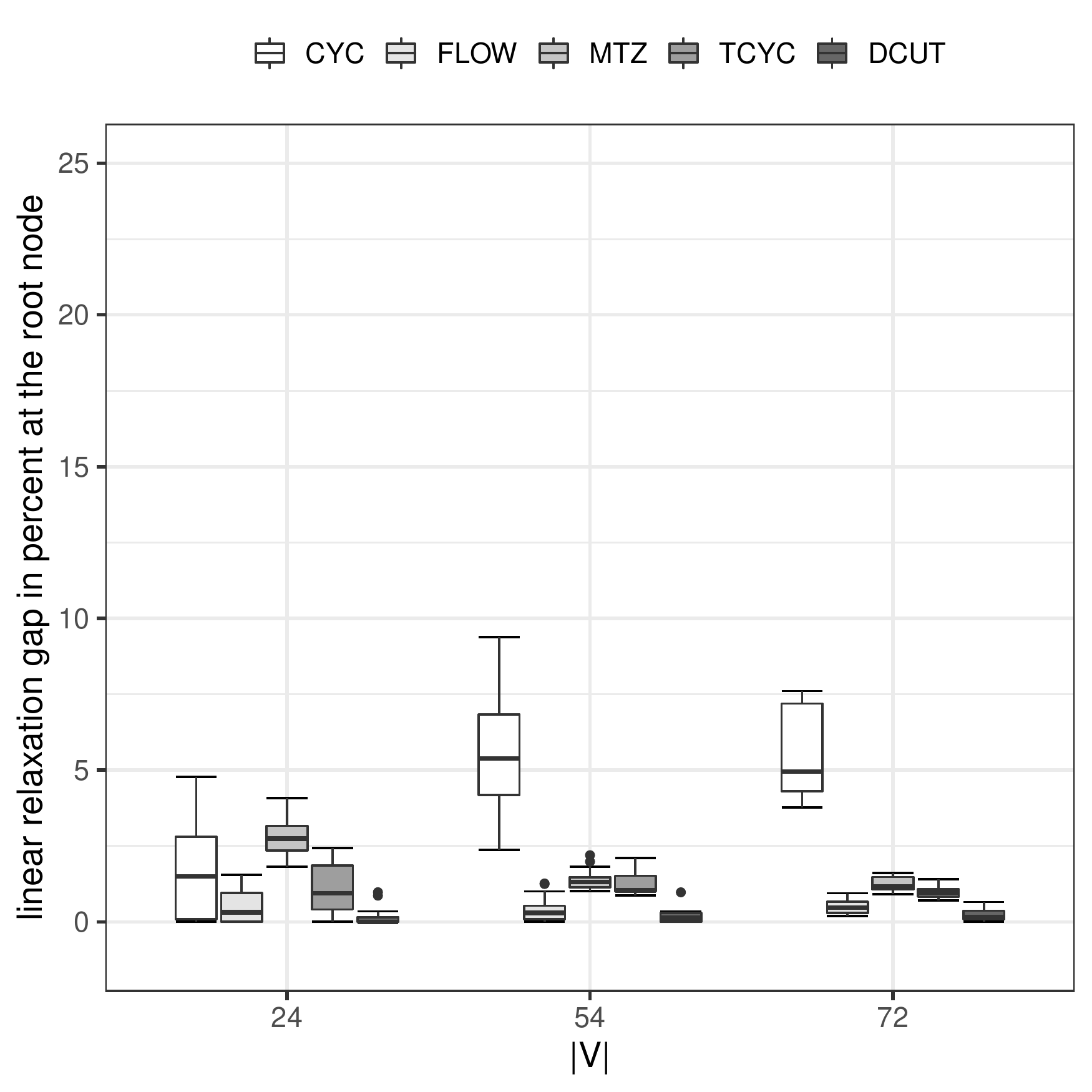}
   \label{fig:glrboxplotb}
 }
  \subfigure[H]{
   \includegraphics[scale =0.45] {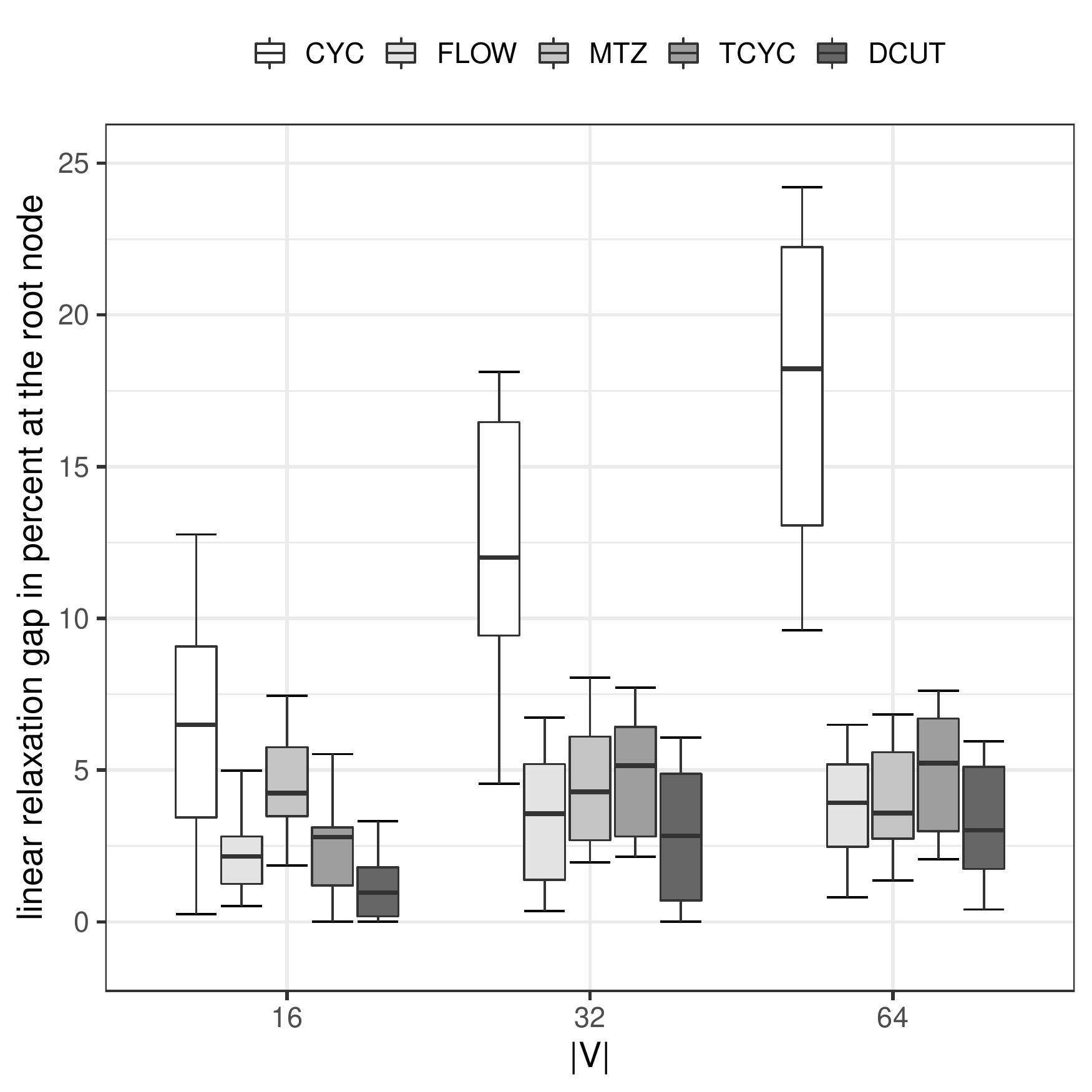}
   \label{fig:glrboxplotc}
 }
  \subfigure[T]{
   \includegraphics[scale =0.45] {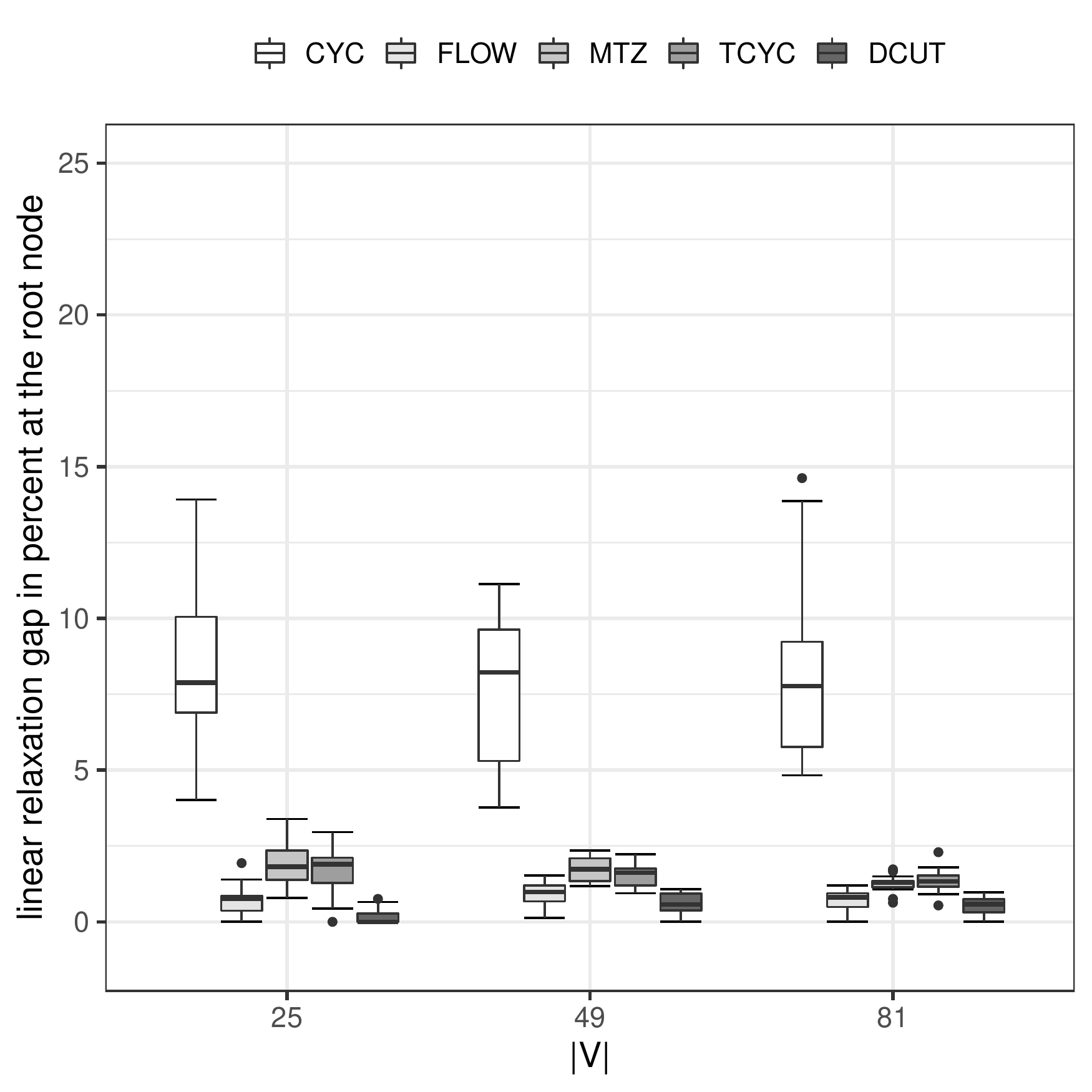}
   \label{fig:glrboxplotd}
 }
\caption{Boxplot summarizing the linear relaxation gaps at the end of the root node (glr) by each formulation for the square grid, nonsquare grid, hypercube, and toroidal instances.}
\label{fig:glrboxplota-d}
\end{figure}

\begin{figure}[H]
\centering
   \includegraphics[width=0.5\textwidth] {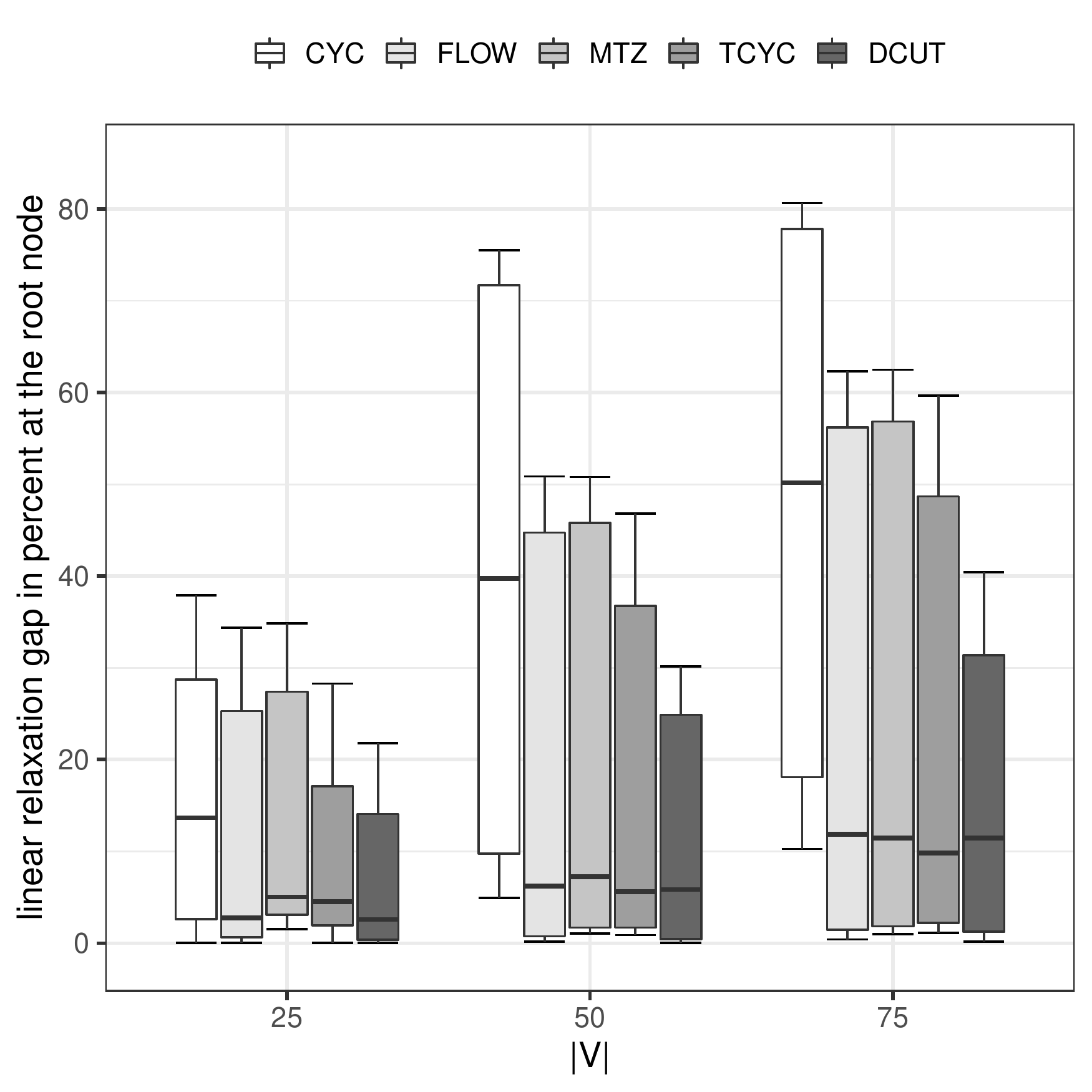}
\caption{Boxplot summarizing the linear relaxation gaps at the end of the root node (glr) by each formulation for the random instances.}
\label{fig:glrboxplote}
\end{figure}

The boxplots in Figures~\ref{fig:glrboxplota-d}-\ref{fig:glrboxplote} summarize the results regarding the linear relaxation gaps at the end of the root node for each of the formulations, separated by the number of vertices.

\section{Extension to the maximum weighted induced tree problem}
\label{sec:adapting}

In this section, we first show that the formulations proposed for problem MWIF in Section~\ref{sec:formulations} can be extended to the maximum weighted induced tree problem (MWIT). After that, we perform computational experiments to illustrate the difference between the optimal solution values of these two problems and provide some examples of structures that imply such differences.

\subsection{Adapting the formulations}

We now show that the tree with cycle elimination, flow-based, MTZ-based and directed cutset formulations can be directly adapted to tackle the maximum weighted induced tree problem. 
This can be achieved with the addition of a single constraint. More specifically, the undirected tree with cycle elimination formulation can ensure that the solution is an induced tree by adding the constraint
\begin{equation}
    \sum_{e \in \delta(s)} x_e \leq 1.
\end{equation}
Similarly, the directed flow-based, MTZ-based and cutset formulations can guarantee that the solution is an induced tree by adding the constraint
\begin{equation}
    \sum_{(s,v) \in A^+ (s)} x_{sv} \leq 1.
\end{equation}

\subsection{Comparing the optimal solutions of MWIF and MWIT}

For the results presented in this section, we limit ourselves to the instances with at most 81 vertices. We remark that these instances could all be solved to optimality within the time limit of 3600 seconds for both problems (MWIF and MWIT).
Denote by $z^*_{MWIF}$ and $z^*_{MWIT}$, correspondingly, the optimal solution values for MWIF and MWIT for a given input graph.

\begin{figure}[H]
    \centering
\begin{tikzpicture}[node distance=2cm]
\tikzstyle{vertex}=[circle, draw,line width=1pt,minimum size = 0.65cm]
\tikzstyle{small}=[circle, draw,line width=1pt]
\tikzstyle{inv}=[circle]
\tikzstyle{edge} = [draw,thick,-,line width=1pt]
\tikzstyle{arcs} = [draw,thick,->,line width=1pt]
	
\node[inv] (s) {};
\node[vertex] [above left of = s] (a) {$a$};
\node[inv] [above of = a, node distance=0.6cm] (l1) {$10$};
\node[vertex] [right of = a] (b) {$b$};
\node[inv] [above of = b, node distance=0.6cm] (l2) {$1$};
\node[vertex] [right of = b] (c) {$c$};
\node[inv] [above of = c, node distance=0.6cm] (l3) {$10$};
\node[vertex] [right of = c] (d) {$d$};
\node[inv] [above of = d, node distance=0.6cm] (l4) {$10$};
\node[vertex] [below of = a] (e) {$e$};
\node[inv] [above left of = e, node distance=0.6cm] (l5) {$10$};
\node[vertex] [right of = e] (f) {$f$};
\node[inv] [above left of = f, node distance=0.6cm] (l6) {$10$};
\node[vertex] [right of = f] (g) {$g$};
\node[inv] [above left of = g, node distance=0.6cm] (l7) {$1$};
\node[vertex] [right of = g] (h) {$h$};
\node[inv] [above left of = h, node distance=0.6cm] (l8) {$10$};
\node[vertex] [below of = e] (i) {$i$};
\node[inv] [below of = i, node distance=0.6cm] (l5) {$10$};
\node[vertex] [right of = i] (j) {$j$};
\node[inv] [below of = j, node distance=0.6cm] (l6) {$1$};
\node[vertex] [right of = j] (k) {$k$};
\node[inv] [below of = k, node distance=0.6cm] (l7) {$10$};
\node[vertex] [right of = k] (l) {$l$};
\node[inv] [below of = l, node distance=0.6cm] (l8) {$10$};

\path[edge] (a) edge [above]node{} (b);
\path[edge] (b) edge [above]node{} (c);
\path[edge] (c) edge [above]node{} (d);

\path[edge] (e) edge [above]node{} (f);
\path[edge] (f) edge [above]node{} (g);
\path[edge] (g) edge [above]node{} (h);

\path[edge] (i) edge [above]node{} (j);
\path[edge] (j) edge [above]node{} (k);
\path[edge] (k) edge [above]node{} (l);

\path[edge] (a) edge [above]node{} (e);
\path[edge] (b) edge [above]node{} (f);
\path[edge] (c) edge [above]node{} (g);
\path[edge] (d) edge [above]node{} (h);

\path[edge] (e) edge [above]node{} (i);
\path[edge] (f) edge [above]node{} (j);
\path[edge] (g) edge [above]node{} (k);
\path[edge] (h) edge [above]node{} (l);

\end{tikzpicture}    
\caption{Example of grid structure forbidding MWIF and MWIT to have the same optimal solution value.}
\label{fig:nqggraph}
\end{figure}
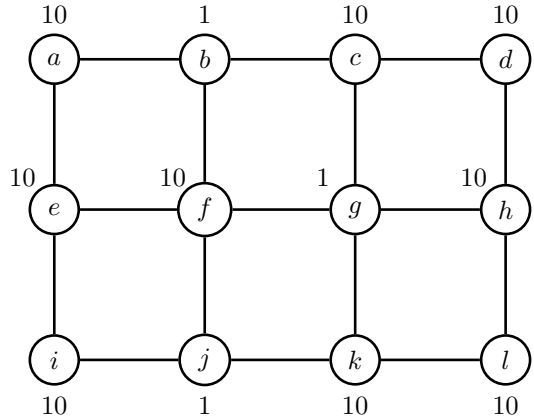

Table~\ref{tab:optimalvaluesmwifmwit} compares the optimal solution values for MWIF and MWIT. In this table, the first two columns identify the graph class and the number of vertices, respectively. The next column gives the total number of instances corresponding to that row. The fourth column (\#diff) provides the number of instances for which $z^*_{MWIF} \neq z^*_{MWIT}$, while the last column (diff (\%)) indicates the average difference for those instances with different optimal solution values for the two problems, calculated as $100\times \frac{z^*_{MWIF}-z^*_{MWIT}}{z^*_{MWIF}}$.
The results indicate that differences in the optimal solutions value occurred in all graph classes, with the exception of the nonsquare grid instances. It is noteworthy that largest values for both the number of differences and the percentual difference occurred in the hypercube instances.
More instances with different optimal values for the two problems were observed for the random graphs than for the other classes (toroidal, grid and nonsquare grid).
For both the toroidal and square grid, differences were only observed for two instances, with considerably low averages.
No difference was observed for the nonsquare grid instances, although there exist such graphs for which the optimal solutions for MWIF and MWIT are different, as the one illustrated later in Figure~\ref{fig:nqggraph}.

\begin{table}[ht]
\centering
\small
\caption{Results comparing the optimal solution values for MWIF and MWIT.}\label{tab:optimalvaluesmwifmwit}
\begin{tabular}{llrrr}
  \hline
Graph class & $|V|$ & Number & \#diff & diff (\%) \\ 
  \hline
Square grid & 25 & 15 & 0 &  \\ 
   & 49 & 15 & 0 &  \\ 
   & 81 & 15 & 2 & 0.3 \vspace{0.2cm} \\ 
  Nonsquare grid & 24 & 15 & 0 &  \\ 
   & 54 & 15 & 0 &  \\ 
   & 72 & 15 & 0 &  \vspace{0.2cm} \\ 
  Hypercube & 16 & 15 & 9 & 4.7 \\ 
   & 32 & 15 & 9 & 1.7 \\ 
   & 64 & 15 & 13 & 1.2\vspace{0.2cm}\\ 
  Toroidal & 25 & 15 & 0 &  \\ 
   & 49 & 15 & 2 & 0.8 \\ 
   & 81 & 15 & 0 &  \vspace{0.2cm} \\ 
  Random & 25 & 45 & 3 & 1.9 \\ 
   & 50 & 45 & 6 & 0.9 \\ 
   & 75 & 45 & 10 & 0.3 \\ 
   \hline
\end{tabular}
\end{table}

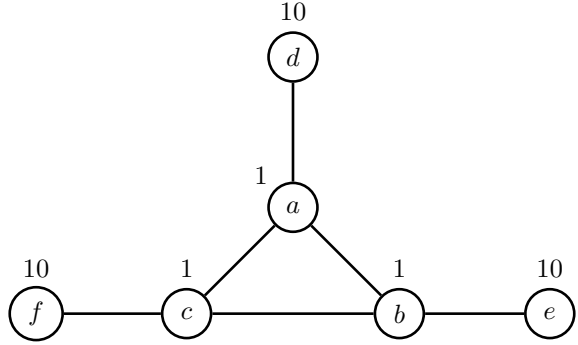
\begin{figure}[H]
    \centering
\begin{tikzpicture}[node distance=2cm]
\tikzstyle{vertex}=[circle, draw,line width=1pt,minimum size = 0.65cm]
\tikzstyle{small}=[circle, draw,line width=1pt]
\tikzstyle{inv}=[circle]
\tikzstyle{edge} = [draw,thick,-,line width=1pt]
\tikzstyle{arcs} = [draw,thick,->,line width=1pt]
	
\node[inv] (s) {};
\node[vertex] [above of = s] (a) {$a$};
\node[inv] [above left of = a, node distance=0.6cm] (l1) {$1$};
\node[vertex] [below right of = a] (b) {$b$};
\node[inv] [above of = b, node distance=0.6cm] (l2) {$1$};
\node[vertex] [below left of = a] (c) {$c$};
\node[inv] [above of = c, node distance=0.6cm] (l3) {$1$};
\node[vertex] [above of = a] (d) {$d$};
\node[inv] [above of = d, node distance=0.6cm] (l4) {$10$};
\node[vertex] [right of = b] (e) {$e$};
\node[inv] [above of = e, node distance=0.6cm] (l5) {$10$};
\node[vertex] [left of = c] (f) {$f$};
\node[inv] [above of = f, node distance=0.6cm] (l6) {$10$};

\path[edge] (a) edge [above]node{} (b);
\path[edge] (b) edge [above]node{} (c);
\path[edge] (c) edge [above]node{} (a);

\path[edge] (a) edge [above]node{} (d);
\path[edge] (b) edge [above]node{} (e);
\path[edge] (c) edge [above]node{} (f);

\end{tikzpicture}    
\caption{Example of clique structure forbidding MWIF and MWIT to have the same optimal solution value.}
\label{fig:3cliquegraph}
\end{figure}

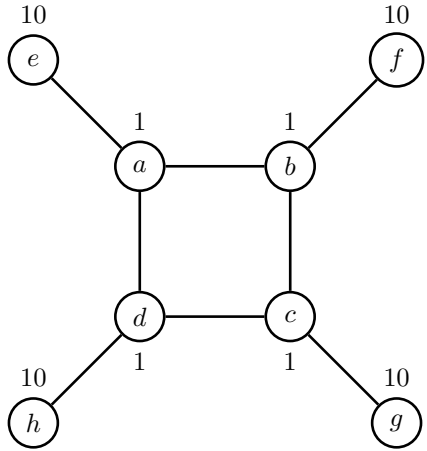
\begin{figure}[H]
    \centering
\begin{tikzpicture}[node distance=2cm]
\tikzstyle{vertex}=[circle, draw,line width=1pt,minimum size = 0.65cm]
\tikzstyle{small}=[circle, draw,line width=1pt]
\tikzstyle{inv}=[circle]
\tikzstyle{edge} = [draw,thick,-,line width=1pt]
\tikzstyle{arcs} = [draw,thick,->,line width=1pt]
	
\node[inv] (s) {};
\node[vertex] [above left of = s] (a) {$a$};
\node[inv] [above of = a, node distance=0.6cm] (l1) {$1$};
\node[vertex] [right of = a] (b) {$b$};
\node[inv] [above of = b, node distance=0.6cm] (l2) {$1$};
\node[vertex] [below of = b] (c) {$c$};
\node[inv] [below of = c, node distance=0.6cm] (l3) {$1$};
\node[vertex] [left of = c] (d) {$d$};
\node[inv] [below of = d, node distance=0.6cm] (l4) {$1$};
\node[vertex] [above left of = a] (e) {$e$};
\node[inv] [above of = e, node distance=0.6cm] (l5) {$10$};
\node[vertex] [above right of = b] (f) {$f$};
\node[inv] [above of = f, node distance=0.6cm] (l6) {$10$};
\node[vertex] [below right of = c] (g) {$g$};
\node[inv] [above of = g, node distance=0.6cm] (l7) {$10$};
\node[vertex] [below left of = d] (h) {$h$};
\node[inv] [above of = h, node distance=0.6cm] (l8) {$10$};

\path[edge] (a) edge [above]node{} (b);
\path[edge] (b) edge [above]node{} (c);
\path[edge] (c) edge [above]node{} (d);
\path[edge] (d) edge [above]node{} (a);

\path[edge] (a) edge [above]node{} (e);
\path[edge] (b) edge [above]node{} (f);
\path[edge] (c) edge [above]node{} (g);
\path[edge] (d) edge [above]node{} (h);

\end{tikzpicture}    
\caption{Example of cycle structure forbidding MWIF and MWIT to have the same optimal solution value.}
\label{fig:cyclegraph}
\end{figure}

\begin{figure}[H]
    \centering
\begin{tikzpicture}[node distance=2cm]
\tikzstyle{vertex}=[circle, draw,line width=1pt,minimum size = 0.65cm]
\tikzstyle{small}=[circle, draw,line width=1pt]
\tikzstyle{inv}=[circle]
\tikzstyle{edge} = [draw,thick,-,line width=1pt]
\tikzstyle{arcs} = [draw,thick,->,line width=1pt]
	
\node[inv] (s) {};
\node[vertex] [above left of = s] (a) {$a$};
\node[inv] [above of = a, node distance=0.6cm] (l5) {$1$};
\node[vertex] [right of = a] (b) {$b$};
\node[inv] [above of = b, node distance=0.6cm] (l6) {$10$};
\node[vertex] [below of = b] (c) {$c$};
\node[inv] [below of = c, node distance=0.6cm] (l2) {$1$};
\node[vertex] [below of = a] (d) {$d$};
\node[inv] [below of = d, node distance=0.6cm] (l3) {$10$};
\node[vertex] [left of = a] (e) {$e$};
\node[inv] [above of = e, node distance=0.6cm] (l1) {$10$};
\node[vertex] [right of = b] (f) {$f$};
\node[inv] [above of = f, node distance=0.6cm] (l4) {$10$};

\path[edge] (a) edge [above]node{} (b);
\path[edge] (b) edge [above]node{} (c);
\path[edge] (c) edge [above]node{} (d);
\path[edge] (d) edge [above]node{} (a);

\path[edge] (a) edge [above]node{} (e);
\path[edge] (d) edge [above]node{} (e);

\path[edge] (b) edge [above]node{} (f);
\path[edge] (c) edge [above]node{} (f);

\end{tikzpicture}    
\caption{Example of combined structure forbidding MWIF and MWIT to have the same optimal solution value.}
\label{fig:combinedgraph}
\end{figure}
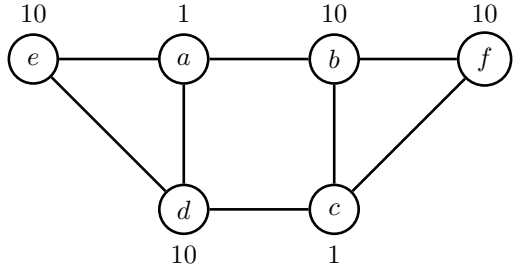

A natural question that arises is related to which types of structures make a difference.
We remark that providing an in-depth study and characterization of structures that forbid achieving the same optimal solution values for MWIF and MWIT is beyond the scope of our work. However, we provide a few examples in Figures~\ref{fig:3cliquegraph}-\ref{fig:nqggraph}.

In the graph depicted in Figure~\ref{fig:3cliquegraph}, an optimal induced forest has weight $z^*_{MWIF}=32$ and one of such forests is composed of two induced trees, one induced by $\{a,c,d,f\}$ with weight 22 and the other induced by $\{e\}$ with weight 10. One can observe that the removal of the clique $\{a,b,c\}$ disconnects the graph in three connected components, implying that the vertices $d$, $e$, and $f$ cannot be together in a solution for MWIT. Thus, a solution for MWIT with $z^*_{MWIT} = z^*_{MWIF}$ cannot be reached. In fact, this happens for any clique $K$ whose removal disconnects the graph in at least three connected components for which reaching the same optimal value for MWIF and MWIF would require having at least three vertices of $K$ to connect all these components.

In the graph illustrated in Figure~\ref{fig:cyclegraph}, an optimal induced forest has weight $z^*_{MWIF}=43$ and one of such forests is composed of two induced trees, one induced by $\{a,b,d,e,f,h\}$ with weight 23 and the other induced by $\{g\}$ with weight 10. Notice that the removal of the cicle $\{a,b,c,d\}$ disconnects the graph in four connected components, implying that the vertices $e$, $f$, $g$, and $h$ cannot be all in the solution to MWIT.
This implies that we cannot obtain a solution for MWIT with $z^*_{MWIT} = z^*_{MWIF}$.
We remark that although this structure seems similar to the one illustrated in Figure~\ref{fig:3cliquegraph}, it relies on the existence of cycles rather than cliques.
 Notice that such prohibition arises for any cycle $C$ whose removal disconnects the graph in at least $|C|$ connected components in a way that achieving $z^*_{MWIF}=z^*_{MWIT}$ would imply using all the vertices of $C$ to connect these components.

Figure \ref{fig:combinedgraph} shows an example in which the way the different cliques and cycles interact with each other imply the impossibility of having the same optimal solution value for MWIF and MWIT. 
The optimal induced forest has weight $z^*_{MWIF}=40$ and is composed of two trees, one induced by $\{b,f\}$ and the other by $\{d,e\}$, each with weight 20. Notice that, given the two cliques of size three, $\{a,d,e\}$ and $\{b,c,f\}$, one can select at most 4 vertices for both MWIF and MWIT. In order to obtain a tree induced by many large weighted vertices in the example, either vertex $a$ or vertex $c$ should be selected, implying that it is not possible to achieve $z^*_{MWIT} = z^*_{MWIF}$.

The graph displayed in Figure~\ref{fig:nqggraph} provides an example whose optimal induced forest has weight $z^*_{MWIF}=90$ and is composed of two induced trees, one induced by $\{a,e,f,i\}$ with weight 40 and the other induced by $\{c,d,h,k,l\}$ with weight 50. Observe that none of the three possible vertices that would be needed to connected these two trees, which are those in $\{b,g,k\}$, could be used without removing one of those in the induced forest.
Therefore, we cannot achieve a solution for MWIT with $z^*_{MWIT} = z^*_{MWIF}$.

\section{Concluding remarks}
\label{sec:concluding}

In this paper, we considered the maximum weighted induced forest problem (MWIF), which given a vertex weighted graph consists of encountering a subset of its vertices with maximum weight inducing a forest. We proposed two new mixed integer programming (MIP) formulations with an exponential number of constraints together with branch-and-cut procedures. 

Extensive computational experiments were performed to compare five MIP formulations for MWIF: two of which are compact, while the three others are based on an exponential number of constraints. More specifically, we compared three existing formulations, namely, a cycle elimination formulation (CYC), a compact flow-based formulation (FLOW), and a compact Miller-Tucker-Zemlin-based formulation (MTZ), and the two new formulations proposed in our work, which are a tree with cycle elimination formulation (TCYC) and a directed cutset formulation (DCUT).

The performed computational experiments have shown that the newly proposed formulations TCYC and DCUT, as well as the compact formulations MTZ and FLOW, recently proposed in~\citeA{MelQueRib21}, clearly outperform CYC, which is based on the existing formulation for the minimum weighted feedback vertex set problem \cite{BruMafTru00}.
Furthermore, the new branch-and-cut approaches have shown to be very effective when tackling the random instances. FLOW, on the other hand, has shown to be a good choice for solving grid, hypercube and toroidal instances. 

We have also shown how certain formulations for MWIF can be easily extended to tackle the maximum weighted induced tree problem (MWIT).
Such extensions allowed us to evaluate the differences in the optimal solution values of MWIF and MWIT when applied to some of the benchmark instances.
The results showed that, considering the benchmark set, the graph class for which the largest percentage of instances presented different optimal values for the two problems was that of the hypercube graphs, followed by random graphs. Besides, the differences between these values were more considerable for this same class of hypercube graphs.
We believe that the proposed approaches can be of great value for graph theoreticians, as they can potentially allow one to verify conjectures.

\section*{Acknowledgments}

Work of Rafael A. Melo was supported by the State of Bahia Research Foundation (FAPESB); and the Brazilian National Council for Scientific and Technological Development (CNPq).
Work of Celso C. Ribeiro was partially supported by CNPq research grants
303958/2015-4 and 425778/2016-9 and by FAPERJ research grant E-26/202.854/2017.
This work was also partially sponsored by Coordena\c{c}\~ao de Aperfeiçoamento de Pessoal de N\'\i vel Superior (CAPES), under Finance Code 001.

\bibliographystyle{apacite}
\bibliography{main}

\end{document}